\documentstyle[12pt]{article}
\input{mssymb.tex}

\catcode`\@=12
\catcode`\~=12

\renewcommand{\thefootnote}{\fnsymbol{footnote}}

\font\goth=eufm10 scaled\magstep1
\def\A{\mbox{\goth A}}
\def\U{\mbox{\goth U}}

\def\ie{\mbox{{\it i.e.} }}
\def\etc{\mbox{{\it etc.} }}
\def\eg{\mbox{{\it e.g.} }}
\def\cf{\mbox{{\it cf.} }}
\def\tr{\mbox{tr}}
\def\id{\mathop{\rm id}}
\def\field#1{{\Bbb #1}}
\def\real{{\field R}}
\def\complex{{\field C}}
\def\R{{\cal R}}
\def\Rhat{\widehat{R}}
\def\matrix#1#2#3{#1 ^{#2}{}_{#3}}
\def\ket#1{\left| #1\right\rangle}
\def\abs#1{\left| #1\right|}
\def\half{{1\over 2}}
\def\inv#1{{1\over #1}}
\def\comm#1#2{\left[#1, #2\right]}
\def\inprod#1#2{\left\langle #1, #2\right\rangle}
\def\trg{{\triangleright}}
\def\tlg{{\triangleleft}}
\def\D#1{{{\cal D}}^#1}
\def\K#1#2{K^{#1}{}_{#2}(\rho)}
\def\M{{\cal V}}
\def\IA{1}
\def\IU{1}

\def\DA{{\Delta_{\A}}}
\def\AD{{{}_{\A}\Delta}}
\def\ad{{\stackrel{\mbox{\scriptsize ad}}{\triangleright}}}
\def\Ad{{\Delta^{\mbox{\scriptsize Ad}}}}
\def\killing#1#2#3{\eta^{(#1)}\left( #2 \otimes #3\right) }
\def\metric#1#2{\eta^{(#1)}_{#2}}
\def\invmet#1#2{\eta^{(#1)\,{#2}}}
\def\trace#1#2{\tr_{#1}\left( u #2 \right)}
\def\k{k}
\def\O#1#2{O_{#1}{}^{#2}}
\def\bigR#1#2{\matrix{\widehat{\real}}{#1}{#2}}
\def\bigA#1#2{\matrix{\field{A}}{#1}{#2}}
\def\bigD#1#2{\matrix{\field{D}}{#1}{#2}}
\def\struc#1#2{f_{#1}{}^{#2}}
\def\pstruc#1#2{f'_{#1}{}^{#2}}
\def\gen#1{T_{#1}}
\def\pgen#1{T'_{#1}}
\def\I#1#2{\tau_{#2}\left( #1 \right)}
\def\pI#1{\tau'_{#1}}
\def\quint#1#2{\left[ #1 \right] _{#2}}
\def\cas#1{Q^{( #1 )}}
\def\pcas{Q'}
\def\adj{ad}
\def\padj{ad'}
\def\fun{fn}
\def\index#1{{\cal I}\left( #1\right)}
\def\g{\mbox{\goth g}}
\def\uea#1#2{U_q(\mbox{\goth #1}( #2 ))}

\begin{document}
\begin{titlepage}
\begin{center}

Centre de Physique Th\'eorique\footnote{Unit\'e Propre de Recherche 7061} -
CNRS - Luminy, Case 907\\
F-13288 Marseille Cedex 9 - France

\vspace{2 cm}

{\large \bf Killing Form on Quasitriangular Hopf Algebras and Quantum Lie
Algebras}

\vspace{1 cm}

{\bf Paul WATTS}\footnote{E-mail: {\it watts@cpt.univ-mrs.fr} - WWW Home
Page: {\it http://cpt.univ-mrs.fr/~watts}}

\vspace{2 cm}

{\bf Abstract}

\end{center}

\noindent The basics of quasitriangular Hopf algebras and quantum Lie
algebras are briefly reviewed, and it is shown that their properties allow
the introduction of a Killing form.  For quantum Lie algebras, this leads
to the definitions of a Killing metric and quadratic casimir.  The specific
case of $\uea{su}{N}$ is examined in detail, where it is shown that many of
the classical results are reproduced, and explicit calculations to
illustrate the conclusions are presented for $\uea{su}{2}$.

\vspace{2 cm}

\noindent MSC-91: 16W30, 17B37, 81R50\\
\noindent PACS-95: 02.20.Sv

\bigskip

\noindent May 1995

\noindent CPT-95/P.3201
%\noindent hep-th/
%\noindent q-alg
\bigskip

\noindent Anonymous FTP or gopher: {\it cpt.univ-mrs.fr}

\end{titlepage}
\setcounter{page}{1}
\renewcommand{\thepage}{\roman{page}}
\tableofcontents
\newpage
\setcounter{page}{1}
\renewcommand{\thepage}{\arabic{page}}

\setcounter{footnote}{0}
\renewcommand{\thefootnote}{\arabic{footnote}}

\section{Introduction}
\setcounter{equation}{0}

In studying an algebra $\g$, one often finds it useful to consider the
Killing metric on $\g$.  From a mathematical perspective, this allows some
insight into some of the algebra's properties.  For instance, $\g$ is a
compact Lie group if the Killing metric associated with the adjoint
representation is positive definite.  As another example, it is often the
case that the Killing metric computed from the trace over a given
representation (``rep'') may say something about its reducibility or
irreducibility.  It also figures prominently in the analysis of certain
infinite-dimensional Lie algebras; if one has a finite-dimensional Lie
algebra $\g$, then the central extension which appears in the commutation
relations of the affine algebra $\widehat{\g}$ depends explicitly upon the
Killing metric of $\g$.

\bigskip

However, from a physicist's point of view, the importance of this object
does not lie so much in its usefulness as a way of studying the abstract
structure of an algebra, but rather as a tool in constructing a physical
theory.  Perhaps the most obvious example of this is when one wants to find
quadratic casimirs; recall that the quadratic casimir of a Lie algebra $\g$
is central within the universal enveloping algebra $U(\g)$.  Therefore,
given an irreducible representation (``irrep'') of $\g$, the matrix
corresponding to the quadratic casimir is proportional to the identity, and
the constant of proportionality labels the irrep.  In constructing a
physical theory for which $\g$ is a symmetry algebra, one (generally)
requires that a particle in the spectrum of our theory lives in an irrep of
$\g$, and therefore the good quantum numbers describing this particle will
include not only the weights of the irrep, but also the value of the
quadratic casimir.  (The most familiar example of this is $\g=${\goth
su}$(2)$, where the quadratic casimir is $\vec{J}^2$, which of course takes
the value $j(j+1)$ in the spin $j$ irrep.)  These casimirs are constructed
(for the cases where $\g$ is semisimple) using the Killing metric, and
therefore the utility of finding such a metric is obvious.

Another very important way in which the Killing metric appears is in the
construction of lagrangeans, in the following fashion: suppose one wishes
to construct a Yang-Mills theory in which $\g$ is the Lie algebra of the
gauge group.  If $\Gamma$ is the gauge field 1-form, and $\rho$ the rep in
which it lives, the action is (up to a multiplicative constant)
\begin{equation}
S=\int\tr_{\rho}(F\wedge\star F)+\ldots ,
\end{equation}
where $F=\mbox{d}\Gamma+\Gamma\wedge\Gamma$ is the field strength 2-form.
If $\Gamma\equiv\rho(\gen{A})\Gamma^A$, where $\{\gen{A}\}$ are the
generators of $\g$, then one sees the appearance of the Killing metric
$\metric{\rho}{AB}:=\tr_{\rho}\left(\gen{A}\gen{B}\right)$ associated with
this rep:
\begin{equation}
S=\int\metric{\rho}{AB}F^A\wedge\star F^B+\ldots
\end{equation}
Thus, knowledge of the Killing metric is vital to the construction of a
Yang-Mills theory, since it shows up explicitly in the gauge field kinetic
energy term.

\bigskip

The purpose of this work is to extend the concept of a Killing metric to
the case where the algebra in question is a quasitriangular Hopf algebra.
As just stated, such a metric is of extreme utility when one considers a
physical theory which has some global or local symmetry, so if one wants to
formulate such a theory where the symmetry algebra is a {\em quantum}
rather than a classical one, and examine its particle content, scattering
amplitudes, {\it etc.}, the Killing metric will be of great importance.

The sections immediately following this Introduction will serve principally
to establish the language and notation which appear throughout this paper,
namely, those of Hopf algebras and quantum Lie algebras, and should be
treated as brief reviews, since the majority the material therein may be
found in more detail elsewhere.  The main thrust of this work starts in
Section \ref{chap-Killing}, where the Killing form on an arbitrary
quasitriangular Hopf algebra is first introduced, and then the Killing
metric is found by specifying to the case of a quantum Lie algebra.
$\uea{su}{N}$, and specifically $\uea{su}{2}$, will then be treated in
order to illustrate the connections with the classical case and provide
concrete examples of the results.  In addition, there is an Appendix
containing useful information about the element $u$ which plays such an
important role in the formulation of the Killing form.

\section{Quasitriangular Hopf Algebras}\label{chap-Hopf}
\setcounter{equation}{0}

Not surprisingly, given the subject of this paper, it will be necessary to
review the basics of the representation theory of Hopf algebras and
quasitriangular Hopf algebras (the latter of which gives rise to quantum
groups).  First, a very brief reminder: a {\em Hopf algebra (HA)} $\U$ is a
unital associative algebra (with unit $\IU$) over a field
$\k$\footnote{Throughout this work, the term ``numerical'' is equivalent to
``$\k$-valued''.} on which there exist the following linear maps: the {\em
coproduct} (or {\em comultiplication}) $\Delta:\U\rightarrow\U\otimes\U$,
the {\em counit} $\epsilon:\U\rightarrow\k$, and the {\em antipode}
$S:\U\rightarrow\U$, which satisfy the conditions
\begin{eqnarray}
(\Delta\otimes\id)\Delta(x)=(\id\otimes\Delta)\Delta(x),&&\Delta(xy)=
\Delta(x)\Delta(y),\nonumber\\
(\epsilon\otimes\id)\Delta(x)=(\id\otimes\epsilon)\Delta(x)=x,&&
\epsilon(xy)=\epsilon(x)\epsilon(y),\nonumber\\
\cdot(S\otimes\id)\Delta(x)=\cdot(\id\otimes S)\Delta(x)=\IU\epsilon(x),&&
\end{eqnarray}
$x,y\in\U$.  Furthermore, if a *-structure is desired, the involution
$\theta:\U\rightarrow\U$ must be an idempotent antilinear map such that
\begin{eqnarray}
\theta(S(\theta(x)))=S^{-1}(x),&\epsilon(\theta(x))=\epsilon(x)^*,&\Delta
(\theta(x))=(\theta\otimes\theta)(\Delta(x)),
\end{eqnarray}
where $^*$ is the conjugation in $\k$.  (These are just the defining
properties; if the reader is interested in learning more about HAs, s/he
may refer to \cite{Sweedler,Abe,Majid1}.)

\bigskip

Before proceeding to the next subsection, it is necessary to introduce some
terminology: by the ``rep of a HA $\U$'', I will always mean a 1-to-1
linear map $\rho$ from $\U$ to some finite-dimensional matrix group over
the field $\k$ such that $\rho(xy) =\rho(x)\rho(y)$, just as in the
classical case.  A ``irrep'' is a nonreducible rep of $\U$, with
reducibility also being defined in the same sense as in the classical case,
\ie the existence of a common nonzero null eigenvector for all matrices in
the rep.  The reduction of a rep $\rho$ to irreps $\{\rho_i\}$ as a direct
sum $\bigoplus_i\rho_i$ then follows accordingly.

Notice that nowhere in these definitions has there been any mention of
anything but the structure of $\U$ as an {\em algebra}; the fact that $\U$
is a {\em Hopf algebra} appears when one considers tensor product reps of
$\U$, whose constructions are defined via the coproduct, namely, if
$\matrix{\rho}{i}{j}$ and $\matrix{\widehat{\rho}}{A}{B}$ are two reps of
$\U$, then the tensor product rep $\rho\otimes\widehat{\rho}$ is defined by
\begin{equation}
\matrix{(\rho\otimes\widehat{\rho})}{iA}{jB}(x):=\matrix{\rho}{i}{j}(x_{(1
)})\matrix{\widehat{\rho}}{A}{B}(x_{(2)}),
\end{equation}
where $\Delta(x)\equiv x_{(1)}\otimes x_{(2)}$ is the very useful Sweedler
notation for the coproduct \cite{Sweedler}.  Therefore, the HA structure
will emerge when, for example, one wants to find the irreps $\{\rho_i\}$
which appear in the decomposition $\rho\otimes\widehat{\rho}=\bigoplus_i
\rho_i$.

\bigskip

One final note for this subsection: when the deformation parameter $q$ is
introduced later, it will always be treated as a formal variable; in other
words, as an arbitrary element of the field $\k$.  There exist situations,
however, where the actual value of $q$ is important, \eg when $\k=
\complex$, the quantity $1+q^2$ vanishes when $q=\pm i$, and this may cause
certain expressions to blow up, or irreps to suddenly become reducible, or
the like.  In this work, however, these situations will not be considered.

\subsection{Representations of Hopf Algebras}\label{chap-Hopf-reps}

Two HAs (both over the same field $\k$) $\U$ and $\A$ are said to be {\em
dually paired} if there exists a nondegenerate inner product
$\inprod{\,}{\,}:\U\otimes\A\rightarrow\k$ such that their respective
units, coproducts, counits, and antipodes satisfy
\begin{eqnarray}
\inprod{xy}{a}&=&\inprod{x\otimes y}{\Delta(a)},\nonumber\\
\inprod{\IU}{a}&=&\epsilon(a),\nonumber\\
\inprod{\Delta (x)}{a\otimes b}&=&\inprod{x}{ab},\nonumber\\
\epsilon(x)&=&\inprod{x}{\IA},\nonumber\\
\inprod{S(x)}{a}&=&\inprod{x}{S(a)},\label{pairing}
\end{eqnarray}
and if a *-structure exists, the involutions satisfy
\begin{equation}
\inprod{\theta(x)}{a}=\inprod{x}{\theta(S(a))}^*,
\end{equation}
where $x,y\in\U$ and $a,b\in\A$.  (It easy to show that all the relevant
consistency relations to ensure that the two algebras are indeed HAs are
satisfied.)

One may also use this dual pairing constructively: if given a HA $\U$,
together with a $N\times N$ matrix rep $\rho:\U\rightarrow M_N(\k)$, it is
possible to define another HA dually paired with $\U$
\cite{Woronowicz3}; this new HA $\A$ is taken to be the one generated by
the $N^2$ elements $\matrix{A}{i}{j}$ defined through
\begin{equation}
\matrix{\rho}{i}{j}(x)\equiv\inprod{x}{\matrix{A}{i}{j}}.
\end{equation}
The faithfulness of the rep guarantees that this inner product is
nondegenerate, and thus that the elements of the matrix $A$ are uniquely
determined; furthermore, the fact that $\rho$ is an algebra map immediately
gives
\begin{equation}
\Delta(\matrix{A}{i}{j})=\matrix{A}{i}{k}\otimes\matrix{A}{k}{j},\,
\epsilon(\matrix{A}{i}{j})=\delta^i_j ,\,S(\matrix{A}{i}{j})=
\matrix{(A^{-1})}{i}{j}.\label{QG-Hopf}
\end{equation}
The multiplication on $\A$, \ie products between the entries of $A$, will
of course depend upon the explicit form of the coproduct on $\U$, through
the third of (\ref{pairing}).

\subsection{Representations of Quasitriangular Hopf Algebras and Quantum
Groups}

A {\em quasitriangular Hopf algebra (QHA)} \cite{Drinfeld} is a HA $\U$
together with an invertible element, the {\em universal R-matrix},
$\R=r_{\alpha}\otimes r^{\alpha}\in\U\otimes\U$ (summation implied), which
satisfies the relations
\begin{eqnarray}
(\Delta\otimes\id)(\R)&=&\R_{13}\R_{23},\nonumber\\
(\id\otimes\Delta)(\R)&=&\R_{12}\R_{23},\nonumber\\
(\sigma\circ\Delta)(x)&=&\R\Delta(x)\R^{-1},\label{univ-R}
\end{eqnarray}
where $\sigma(x\otimes y)=y\otimes x$, and the subscripts on $\R$ are a
shorthand for
\begin{eqnarray}
\R_{12}&=&r_{\alpha}\otimes r^{\alpha}\otimes\IU,\nonumber\\
\R_{13}& =&r_{\alpha}\otimes\IU\otimes r^{\alpha},\nonumber\\
\R_{23}&=&\IU\otimes r_{\alpha}\otimes r^{\alpha}.
\end{eqnarray}
As a consequence of (\ref{univ-R}), $\R$ satisfies the {\em quantum
Yang-Baxter equation (QYBE)}
\begin{equation}
\R_{12}\R_{13}\R_{23}=\R_{23}\R_{13}\R_{12}.\label{QYBE}
\end{equation}

Of vital importance to the eventual introduction of the Killing form on a
QHA $\U$ will be the element $u$, defined by
\begin{equation}
u:=\cdot (S\otimes\id)(\R_{21})=S(r^{\alpha})r_{\alpha}.\label{u-def}
\end{equation}
$u$ is invertible, with inverse $u^{-1}=r^{\alpha}S^2 (r_{\alpha})$, and
generates the square of the antipode on $\U$ by conjugation, \ie
\begin{equation}
S^2(x)=uxu^{-1}.\label{uxu}
\end{equation}
A consequence of this property is that the element $c:=uS(u)$ is central in
$\U$.

\bigskip

Given a rep $\rho$, the commutation relations between the entries of $A$
may be explicitly determined by using the last of the properties of the
universal R-matrix from (\ref{univ-R}): if
\begin{equation}
\matrix{R}{ij}{k\ell}:=\inprod{\R}{\matrix{A}{i}{k}\otimes
\matrix{A}{j}{\ell}}
\end{equation}
is the $N^2\times N^2$ dimensional {\em numerical R-matrix of} $\A$, one
finds the ``RAA equation'' \cite{RTF}
\begin{equation}
RA_1A_2=A_2A_1R,\label{RAA}
\end{equation}
where the matrix indices have been replaced by the subscripts in an obvious
way.  It is immediate that the QYBE has a numerical counterpart, usually
simply referred to as the {\em Yang-Baxter equation (YBE)}:
\begin{equation}
R_{12}R_{13}R_{23}=R_{23}R_{13}R_{12}\label{YBE}
\end{equation}
(which takes the form $\Rhat_{12}\Rhat_{23}\Rhat_{12}=\Rhat_{23}\Rhat_{12}
\Rhat_{23}$ when one uses the permutation matrix $\matrix{P}{ij}{k\ell}:=
\delta^i_{\ell}\delta^j_k$ to define $\matrix{\Rhat}{ij}{k\ell}:=
\matrix{(PR)}{ij}{k\ell}\equiv\matrix{R}{ji}{k\ell}$).

This leads to the following definition: a HA $\A$ which is dually paired
with a QHA $\U$ by means of a representation $\rho$ in this manner is
called a a {\em matrix pseudogroup} or, more commonly nowadays, a {\em
quantum group (QG)} \cite{Drinfeld}.

\bigskip

(\ref{YBE}) was obtained from (\ref{QYBE}) by applying the given rep to all
three spaces of $\U\otimes\U\otimes\U$, \eg
\begin{equation}
(\rho\otimes\rho\otimes\rho)(\R_{12}\R_{13}\R_{23})=R_{12}R_{13}R_{23}.
\end{equation}
It is also useful to consider the cases where one takes the rep in only one
or two of the tensor product spaces.  To see this, define the matrices
$L^{\pm}\in M_N(\U)$ by
\begin{eqnarray}
L^+&:=&(\id\otimes\rho)(\R)\equiv r_{\alpha}\rho(r^{\alpha}),\nonumber\\
L^-&:=&(\rho\otimes\id)(\R^{-1})\equiv\rho(S(r_{\alpha}))r^{\alpha}.
\label{L-defs}
\end{eqnarray}
{}From the properties of $\R$, one then finds that
\begin{eqnarray}
\Delta(L^{\pm})=L^{\pm}\dot{\otimes}L^{\pm},&&\epsilon(L^{\pm})=I,
\nonumber\\
S(L^+)=(L^+)^{-1}=(\id\otimes\rho)(\R^{-1}),&&S(L^-)=(L^-)^{-1}=(\rho\otimes
\id)(\R)\label{L-hopf}
\end{eqnarray}
($\matrix{(M\dot{\otimes}N)}{i}{j}:=\matrix{M}{i}{k}\otimes
\matrix{N}{k}{j}$).  Now, one can apply the rep to any two of the three
spaces involved in (\ref{QYBE}) to find
\begin{eqnarray}
L^{\pm}_1L^{\pm}_2R=RL^{\pm}_2L^{\pm}_1,&&L^-_2L^+_1R=RL^+_2L^-_1.
\label{LL1}
\end{eqnarray}
The numerical matrices for $L^{\pm}$ also follow:
\begin{eqnarray}
\matrix{\rho}{i}{j}(\matrix{(L^+)}{k}{\ell})&=&\matrix{R}{ik}{j\ell},
\nonumber\\
\matrix{\rho}{i}{j}(\matrix{(L^-)}{k}{\ell})&=&\matrix{(R_{21}^{-1})
}{ik}{j\ell}.
\end{eqnarray}
(As the reader shall see, these will figure very prominently in the
construction of QLAs.)

\section{Actions and Coactions}\label{chap-coactions}
\setcounter{equation}{0}

In the classical case, a Killing form for an algebra $\g$ is defined so
that it is invariant under some action of $\g$ on itself (\eg the adjoint
action in the case of a Lie algebra); therefore, since the goal is to
extend the definition to include QHAs, one must discuss the ways in which a
QHA can act on itself before proceeding further in this direction.  This
requires introducing the concepts of actions and coactions.  (A fuller
discussion of the ideas presented here may be found in \cite{SWZ}.)

\subsection{Actions}\label{chap-coactions-actions}

Suppose $\U$ is a HA and $\M$ a unital associative algebra (both over the
same field $\k$); a {\em (left) action} or {\em generalized (left)
derivation} of $\U$ on $\M$ is a linear map $\trg:\U\otimes\M\rightarrow\M$
satisfying
\begin{eqnarray}
(xy)\trg v=x\trg (y\trg v),&&1\trg v=v,\nonumber\\
x\trg(vw)=(x_{(1)}\trg v)(x_{(2)}\trg w),&&x\trg 1=1\epsilon(x).
\end{eqnarray}
for all $x,y\in\U$ and $v,w\in\M$.  (Note that this is merely another way
of saying that $\trg$ defines a linear rep of $\U$ with right module $\M$.
A right action $\tlg$ of $\U$ on $\M$ can be defined similarly.)  $\U$ may
therefore be interpreted as an algebra of differential operators acting on
elements of $\M$ from the left, and, as such, may be thought of as
providing a means for generalizing infinitesimal transformations.

There are two extremely important examples of such generalized derivations,
both of which will be relevant for this work:
\begin{itemize}
\item The {\em (left) adjoint action} of a Hopf algebra $\U$ on itself is
defined as the linear map $\ad :\U\otimes\U\rightarrow\U$ given by
\begin{equation}
x\ad y:=x_{(1)}yS(x_{(2)}).\label{ad-act}
\end{equation}
In an undeformed Lie algebra, where $\Delta(x)=x\otimes\IU+\IU\otimes x$
and $S(x)=-x$, this reduces to the usual commutator $x\ad y\equiv
\comm{x}{y}$.
\item If $\U$ and $\A$ are two dually paired Hopf algebras, one may
define the (left) action of $\U$ on $\A$ as
\begin{equation}
x\trg a:=a_{(1)}\inprod{x}{a_{(2)}}.\label{diffact}
\end{equation}
As stated above, this allows the interpretation of $\U$ as an algebra of
differential operators which act on elements of $\A$ (``functions'').
\end{itemize}

\subsection{Coactions}\label{chap-coactions-coacts}

$\A$ is a HA and $\M$ is again a unital associative algebra; a {\em
(right) coaction} of $\A$ on $\M$ is a linear map $\DA:\M\rightarrow\M
\otimes\A$ satisfying
\begin{eqnarray}
(\DA\otimes\id)\DA(v)=(\id\otimes\Delta)\DA(v),&&(\id\otimes\epsilon)\DA(v)
=v,
\end{eqnarray}
for all $v\in\M$.  (Many of the following equations will be written using
the ``Sweedleresque'' notation $\DA(v)= v^{(1)}\otimes v^{(2)'}$,
where the unprimed elements live in $\M$, the primed elements in $\A$.)  A
left coaction $\AD(v)=v^{(1)'}\otimes v^{(2)}$ may be defined
similarly.  If $\A_1$ and $\A_2$ are two HAs which coact on $\M$ from the
left and from the right respectively, it will in general be assumed that
they commute with one another, \ie
\begin{equation}
(_{\A_1}\Delta\otimes\id)\Delta_{\A_2}(v)=(\id\otimes\Delta_{\A_2})_{\A_1}
\Delta(v).
\end{equation}
Finally, if $v$ has the property that $\DA(v)=v\otimes 1$, it is said to be
{\em right-invariant}.

One coaction which will appear later in this work is the {\em (right)
adjoint coaction} $\Ad:\A\rightarrow\A\otimes\A$ of $\A$ on itself, defined
via
\begin{equation}
\Ad (a):=a_{(2)}\otimes S(a_{(1)})a_{(3)}.\label{adj-coact}
\end{equation}

\bigskip

A comment on terminology: as the reader may have guessed, the reason for
the term ``coaction'' is because of duality.  If $\U$ and $\A$ are dually
paired HAs, and $\trg$ is a (left) action of $\U$ on some unital
associative algebra $\M$, then there is a natural way to induce a (right)
coaction of $\A$ on $\M$ via
\begin{equation}
v^{(1)}\inprod{x}{v^{(2)'}}=x\trg v.
\end{equation}
The interpretation of this particular coaction is straightforward: it is
the generalization of a {\em finite} transformation of an element of $\M$,
as opposed to the {\em infinitesimal} transformation provided by the
action.

Of particular interest here is the case where $\M=\U$, and the action in
question is the adjoint action $\ad$.  Then the coaction $\DA:\U
\rightarrow\U\otimes\A$ uniquely defined by this pairing, written as $y\mapsto
y^{(1)}\otimes y^{(2)'}$, is given by
\begin{equation}
x\ad y=y^{(1)}\inprod{x}{y^{(2)'}}.\label{coact}
\end{equation}
Notice that this coaction is an algebra map, \ie $\DA(xy)=\DA(x)\DA(y)$.
This will not in general be the case, as the reader can see from the
adjoint coaction $\Ad$ above, which is {\em not} an algebra map from $\A$
to $\A\otimes\A$.

\section{Quantum Lie Algebras}\label{chap-QLA}
\setcounter{equation}{0}

\subsection{Basics of Quantum Lie Algebras}\label{chap-QLA-basics}

Let $\U$ be a HA; it is also a {\em quantum Lie algebra (QLA)} iff there
exists a finite-dimensional subspace $\g\subset\U$ (dim$\,\g=n$) which has
the following properties\footnote{This set of conditions is in effect the
dually paired version of the one required for the existence of a general
differential calculus \cite{Woronowicz1}.}:
\newcounter{QL}
\begin{list}%
{(\Roman{QL})}{\usecounter{QL}
\setlength{\rightmargin}{\leftmargin}}
\item as a vector space, $\U\equiv U_q(\g)$, \ie the universal enveloping
algebra (UEA) of $\g$, modulo commutation relations;
\item the adjoint action $\ad$ closes on $\g$, \ie $\U\ad\g\subseteq\g$;
\item $\Delta(x)\in(\g\otimes\IU)\oplus(\U\otimes\g)$ for all $x\in\g$;
\item for all $x\in\g$, $\epsilon(x)=0$.
\end{list}
The $q$ subscript in (I) simply indicates that the commutation relations
may be deformed relative to the classical case.  (In QG language, $q$ is
the parameter describing the degree of deviation of the QHA from the
undeformed case, the latter of which corresponds to $q=1$.)

Let $\{\gen{A}|A=1,\ldots,n\}$ be a linearly independent basis for $\g$
\cite{Bernard}; the above properties require that
\begin{eqnarray}
\Delta(\gen{A})=\gen{A}\otimes\IU+\O{A}{B}\otimes\gen{B},&\epsilon(\gen{A})
=0,&S(\gen{A})=-S(\O{A}{B})\gen{B},\label{cop-chi}
\end{eqnarray}
where the $n^2$ elements $\O{A}{B}\in\U$ must satisfy
\begin{eqnarray}
\Delta(\O{A}{B})=\O{A}{C}\otimes\O{C}{B},&\epsilon(\O{A}{B})=\delta_A^B,&
S(\O{A}{B})=(O^{-1})_A{}^{B}.\label{O-hopf}
\end{eqnarray}
Condition (II), the closure of $\g$ under the adjoint action of $\U$,
allows the definition of the $\k$-numbers $\bigR{AB}{CD}$ and
$\struc{AB}{C}$ via
\begin{eqnarray}
\gen{A}\ad\gen{B}:=\struc{AB}{C}\gen{C},&&\O{C}{B}\ad\gen{D}:=\bigR{AB}{CD}
\gen{A}\label{ad-chi}
\end{eqnarray}
(from which one can also show $S^{-1}(\O{D}{A})\ad\gen{C}=
\matrix{(\bigR{-1}{})}{AB}{CD}\gen{B}$).  $\bigR{}{}$ is referred to as the
R-matrix of $\g$, and the $f$s are, as in the classical case, just the
structure constants.  Now, note that the definition of the adjoint action
(\ref{ad-act}) implies the identity
\begin{equation}
xy\equiv (x_{(1)}\ad y)x_{(2)};\label{adj-mult}
\end{equation}
therefore, using (\ref{cop-chi}--\ref{ad-chi}), one can show that
\begin{eqnarray}
\gen{A}\gen{B}-\bigR{CD}{AB}\gen{C}\gen{D}&=&\struc{AB}{C}\gen{C},
\nonumber\\
\bigR{EF}{AB}\O{E}{C}\O{F}{D}&=&\O{A}{E}\O{B}{F}\bigR{CD}{EF},\nonumber\\
\gen{A}\O{B}{C}-\bigR{DE}{AB}\O{D}{C}\gen{E}&=&\struc{AB}{D}\O{D}{C}-
\O{A}{D}\O{B}{E}\struc{DE}{C},\nonumber\\
\O{A}{B}\gen{C}&=&\bigR{DE}{AC}\gen{D}\O{E}{B}.\label{Lie-comm}
\end{eqnarray}
The first of these is the deformation of the classical commutation relation
between the generators; the others are due to the difference between the
coproducts of the deformed algebra and the classical versions
$\Delta(\gen{A})=\gen{A}\otimes\IU +\IU\otimes\gen{A}$.

The self-consistency of these relations also requires
\begin{equation}
\bigR{}{12}\bigR{}{23}\bigR{}{12}=\bigR{}{23}\bigR{}{12}\bigR{}{23}.
\label{YBE-QLA}
\end{equation}
(Despite the appearance of a very Yang-Baxterlike equation here, the
numerical matrix $\bigR{}{}$ is {\em not} the rep of a universal R-matrix,
which is good, because so far no mention has been made of any
quasitriangularity of $\U$.)

\subsection{The Adjoint Representation}

The closure of $\g$ under the adjoint action defines the {\em adjoint rep}
$\adj$ of $\U$ in the usual way, \ie
\begin{equation}
x\ad\gen{A}=\matrix{\adj}{B}{A}(x)\gen{B}.
\end{equation}
This motivates the introduction of an $n\times n$ matrix $\bigA{A}{B}$ with
entries in the dually paired HA $\A$, defined through
\begin{equation}
\inprod{x}{\bigA{A}{B}}:=\matrix{\adj}{A}{B}(x).
\end{equation}
It follows that
\begin{eqnarray}
\struc{AB}{C}\equiv\inprod{\gen{A}}{\bigA{C}{B}},&&\bigR{AB}{CD}\equiv
\inprod{\O{C}{B}}{\bigA{A}{D}}.
\end{eqnarray}
Furthermore, by using (\ref{coact}), the coaction of $\A$ on $\U$ is
\begin{equation}
\DA(\gen{A})=\gen{B}\otimes\bigA{B}{A}.\label{gens-coact}
\end{equation}
The coproduct, counit, antipode, and commutation relations on $\A$ have
rather familiar forms:
\begin{eqnarray}
\Delta(\bigA{A}{B})=\bigA{A}{C}\otimes\bigA{C}{B},&&\epsilon(\bigA{A}{B})=
\delta^A_B,\nonumber\\
S(\bigA{A}{B})=\matrix{(\bigA{-1}{})}{A}{B},&&\bigR{}{}\bigA{}{1}\bigA{}{2}
=\bigA{}{1}\bigA{}{2}\bigR{}{},\nonumber\\
\struc{AB}{D}\bigA{C}{D}=\bigA{D}{A}\bigA{E}{B}\struc{DE}{C}.&&
\end{eqnarray}

The above properties of $\bigA{}{}$ imply several numerical relations among
the R-matrix and structure constants: for instance, if one takes the inner
product of $\bigA{M}{N}$ and the first of (\ref{Lie-comm}), the deformed
version of the Jacobi identity appears:
\begin{equation}
\struc{AL}{M}\struc{BN}{L}-\bigR{CD}{AB}\struc{CL}{M}\struc{DN}{L}
=\struc{AB}{C}\struc{CN}{M}.
\end{equation}
Repeating this for the others just recovers (\ref{YBE-QLA}), as well as
\begin{eqnarray}
\bigR{DC}{BN}\struc{AD}{M}-\bigR{DE}{AB}\bigR{MC}{DF}\struc{EN}{F}&=&
\bigR{MC}{DN}\struc{AB}{D}-\bigR{DF}{BN}\bigR{ME}{AD}\struc{EF}{C},
\nonumber\\
\bigR{MB}{AD}\struc{CN}{D}&=&\bigR{DE}{AC}\bigR{FB}{EN}\struc{DF}{M}.
\end{eqnarray}

\subsection{Quasitriangular Quantum Lie Algebras}\label{QQLA}

Given any QHA $\U$ with rep $\rho$, one can immediately obtain a QLA in the
following way \cite{Jurco}: using the matrices $L^{\pm}$, define the matrix
$Y\in M_N(\U)$ \cite{RS,Leipzig} by
\begin{equation}
Y:=L^+S(L^-)\equiv(\rho\otimes\id)(\R_{21}\R)\equiv\rho(r^{\alpha}r_{
\beta})r_{\alpha}r^{\beta}.\label{Y-def}
\end{equation}
As a consequence of (\ref{L-hopf}) and (\ref{LL1}), this matrix satisfies
\begin{eqnarray}
\Delta(\matrix{Y}{i}{j})=\matrix{(L^+)}{i}{k}S(\matrix{(L^-)}{\ell}{j})
\otimes\matrix{Y}{k}{\ell},&&\epsilon(\matrix{Y}{i}{j})=\delta^i_j,
\nonumber\\
S(\matrix{Y}{i}{j})=S^2(\matrix{(L^-)}{k}{j})S(\matrix{(L^+)}{i}{k}),&&
R_{21}Y_1R Y_2=Y_2R_{21}Y_1R.
\end{eqnarray}
To find the coaction of the dually paired HA $\A$ on $Y$, there is the
following trick \cite{SWZ}: for $a\in\A$, define $\Upsilon_a\in\U$ as
\begin{equation}
\Upsilon_a:=\inprod{\R_{21}\R}{a\otimes\id}
\end{equation}
Now, note that for $x\in\U$, one may use the fact that $\R_{21}\R$ is
central in $\Delta(\U)$ to show that
\begin{equation}
x\ad\Upsilon_a=\inprod{x}{S(a_{(1)})a_{(3)}}\Upsilon_{a_{(2)}}.
\label{closure}
\end{equation}
Thus, from (\ref{coact}),
\begin{equation}
\DA(\Upsilon_a)=\Upsilon_{a_{(2)}}\otimes S(a_{(1)})a_{(3)}.
\end{equation}
(Note the appearance of the adjoint coaction (\ref{adj-coact}) in this
equation.)  Since $Y\equiv\Upsilon_A$, where $A$ is the QG matrix, it
immediately follows that
\begin{equation}
\DA(\matrix{Y}{i}{j})=\matrix{Y}{k}{\ell}\otimes S(\matrix{A}{i}{k})
\matrix{A}{\ell}{j}.\label{Y-coact}
\end{equation}
Furthermore, (\ref{closure}) indicates that the adjoint action of $\U$ on
any element in the subspace $\{\Upsilon_a|a\in\A \}$ returns another
element of the same subspace.  In particular,
\begin{equation}
x\ad\matrix{Y}{i}{j}=\inprod{x}{S(\matrix{A}{i}{k})\matrix{A}{\ell}{j}}
\matrix{Y}{k}{\ell},\label{closure-qla}
\end{equation}
which is simply a linear combination of the entries of $Y$.

\bigskip

Notice that in the classical limit, $\R\rightarrow\IU\otimes\IU$, so
$Y\rightarrow I\IU$; therefore, define the matrix $X$ by
\begin{equation}
X:=\frac{I\IU -Y}{\lambda},\label{X-def}
\end{equation}
where $\lambda$ is a parameter which vanishes in the classical limit, and
describes the deviation of $Y$ from unity in the deformed case.  Using the
aforementioned properties of $Y$, one immediately finds
\begin{eqnarray}
\Delta(\matrix{X}{i}{j})=\matrix{X}{i}{j}\otimes\IU+\matrix{(L^+)}{i}{k}S(
\matrix{(L^-)}{\ell}{j})\otimes\matrix{X}{k}{\ell},&&\epsilon (X)=0,
\nonumber\\
S(\matrix{X}{i}{j})=-S^2 (\matrix{(L^-)}{\ell}{j})S(\matrix{(L^+)}{i}{k})
\matrix{X}{k}{\ell},&&\nonumber\\
R_{21}X_1RX_2-X_2R_{21}X_1R=\inv{\lambda}(R_{21}RX_2-X_2R_{21}R),&&
\label{X-X}
\end{eqnarray}
and $\A$ coacts on $X$ as\footnote{Although this work is concerned
primarily with the {\em right} coaction of the QG $\A$ on $\U$, there is
also a {\em left} coaction $x\mapsto\IA\otimes x$, \ie $\U$ is
left-invariant.  The construction presently being discussed reproduces
this, and motivates calling the resulting QLA ``bicovariant''
\cite{SWZ,Woronowicz1}.}
\begin{equation}
\DA(\matrix{X}{i}{j})=\matrix{X}{k}{\ell}\otimes S(\matrix{A}{i}{k})
\matrix{A}{\ell}{j}.\label{X-coacts}
\end{equation}
Therefore, as was the case with $Y$, the adjoint action of $x\in\U$ on an
entry of $X$ returns a linear combination of entries of $X$, so if $\g$ is
defined as the subspace of $\U$ spanned by the entries of $X$ over $\k$,
the UEA $U_q(\g)$ satisfies all four criteria needed for a QLA.  To connect
to the contents of Section \ref{chap-QLA-basics}, one first computes the
adjoint action of one entry of $X$ on another:
\begin{equation}
\matrix{X}{i}{j}\ad\matrix{X}{k}{\ell}=\inv{\lambda}\left[\delta^i_j
\matrix{X}{k}{\ell}-\matrix{(R_{21}^{-1}X_2R_{21}R)}{im}{n\ell}
\matrix{\widetilde{R}}{nk}{jm}\right] \label{X-ad-X}
\end{equation}
(see (\ref{tilde}) for the definition of $\widetilde{R}$).  Now, comparison
of (\ref{X-ad-X}) with (\ref{Lie-comm}) immediately leads to the
identifications
\begin{eqnarray}
\gen{(ij)}\equiv\matrix{X}{i}{j},&\O{(ij)}{(k\ell )}\equiv
\matrix{(L^+)}{i}{k}S(\matrix{(L^- )}{\ell}{j}),&\bigA{(ij)}{(k\ell)}
\equiv S(\matrix{A}{k}{i})\matrix{A}{j}{\ell},\label{adjoint-rep}
\end{eqnarray}
and
\begin{eqnarray}
\bigR{(ab)(cd)}{(ij)(k\ell )}&=&\matrix{\widetilde{R}}{mk}{jn}
\matrix{\Rhat}{sd}{m\ell}\matrix{(\Rhat^{-1})}{ni}{ra}
\matrix{\Rhat}{rb}{sc},\nonumber\\
\struc{(ij)(k\ell )}{(rs)}&=&\inv{\lambda}\left[ \delta^i_j\delta^k_r
\delta^s_{\ell}-\matrix{\widetilde{R}}{mk}{jn}\matrix{(\Rhat^{-1})
}{ni}{tr}\matrix{(\Rhat^2)}{ts}{m\ell}\right] .\label{struc-R}
\end{eqnarray}
There are three important sum relations which follow from these last two
identifications:
\begin{eqnarray}
\struc{AB}{C}I_C=0,&\bigR{CD}{AB}I_C=\delta^D_AI_B,&\bigR{CD}{AB}I_D=
\delta^C_BI_A-\lambda\struc{AB}{C},\label{QLA-I}
\end{eqnarray}
where, for the sake of simplicity, the quantity $I_{(ij)}:=\delta^i_j$ has
been introduced.  There is also the additional relation $I_B\bigA{B}{A}=
I_A\IA$.

(To explain an earlier parenthetical comment, I should point out that the
numerical R-\-mat\-rix in this rep, \ie $\matrix{\real}{AB}{CD}
:=\inprod{\R}{\bigA{A}{C}\otimes\bigA{B}{D}}$, is
\begin{equation}
\matrix{\real}{(ab)(cd)}{(ij)(k\ell )}=\matrix{\widetilde{R}}{mk}{jn}
\matrix{\Rhat}{sb}{m\ell}\matrix{\Rhat}{ni}{rc}
\matrix{(\Rhat^{-1})}{rd}{sa},
\end{equation}
which is {\em not} equal to $\bigR{(cd)(ab)}{(ij)(k\ell )}$.)

\section{The Killing Form}\label{chap-Killing}
\setcounter{equation}{0}

\subsection{The Killing Form for a Quasitriangular Hopf Algebra}

Recall that for a Lie algebra $\g$ over a field $\k$, a Killing form $\eta$
is defined to be a symmetric invariant linear map from $\g\otimes\g$ to
$\k$, \ie $\eta(x\otimes y)=\eta(y\otimes x)$ and $\eta(\comm{z}{x}\otimes
y)+\eta(x\otimes\comm{z}{y})=0$ (\cf \cite{Jac}).  (Such a form always
exists for finite-dimensional Lie algebras: one needs only take the trace
of two elements in the adjoint rep of $\g$.)  Now, with all the formalism
developed so far, it is finally possible to define an object on a QHA which
generalizes these required properties, in the following way:

\bigskip

\noindent {\bf Definition} Let $\U$ be a QHA with universal R-matrix $\R$
and rep $\rho$.  The linear map $\eta^{(\rho)}:\U\otimes\U\rightarrow\k$,
the {\em Killing form associated with} $\rho$, is given by
\begin{equation}
\killing{\rho}{x}{y}:=\trace{\rho}{xy}
\end{equation}
for $x,y\in\U$\footnote{This differs from the definition in \cite{ixtapa}
by an overall antipode, \ie $\killing{\rho}{x}{y}:=\tr_{\rho}\left(
S(uxy)\right)$, which corresponds to taking the trace over the
contragredient rep.}.
\vskip .5cm
$\eta^{(\rho)}$ defined in this way has the following properties:
\begin{eqnarray}
\killing{\rho}{y}{x}&=&\killing{\rho}{x}{S^2(y)},\nonumber\\
\killing{\rho}{(z_{(1)}\ad x)}{(z_{(2)}\ad y)}&=&
\killing{\rho}{x}{y}\epsilon(z),\label{inv-inf}
\end{eqnarray}
for all $x,y,z\in\U$.  The first of these is immediate from the properties
of $u$ and the cyclicity of the trace:
\begin{eqnarray}
\killing{\rho}{y}{x}&=&\trace{\rho}{yx}\nonumber\\
&=&\tr_{\rho}\left( (uyu^{-1})ux\right) \nonumber\\
&=&\trace{\rho}{x(uyu^{-1})}\nonumber\\
&=&\killing{\rho}{x}{S^2(y)}.
\end{eqnarray}
The second is obtained by noting that
\begin{eqnarray}
\trace{\rho}{(z\ad x)}&=&\trace{\rho}{z_{(1)}xS(z_{(2)})}\nonumber\\
&=&\killing{\rho}{z_{(1)}}{xS(z_{(2)})}\nonumber\\
&=&\killing{\rho}{xS(z_{(2)})}{S^2(z_{(1)})}\nonumber\\
&=&\trace{\rho}{xS(S(z_{(1)})z_{(2)})}\nonumber\\
&=&\trace{\rho}{x}\epsilon(z),
\end{eqnarray}
so that by replacing $x$ by $xy$ above, the result follows.

The first of (\ref{inv-inf}) is obviously a generalization of the
``symmetry'' of the Killing form, since $S^2=\id$ in the classical case;
since the adjoint action $\ad$ is the HA generalization of the classical
commutator, the second expresses the desired invariance property of
$\eta^{(\rho)}$.  In fact, this identification may be made more obvious by
realizing that the latter may be rewritten as
\begin{equation}
\killing{\rho}{(z\ad x)}{y}+\killing{\rho}{x}{((-S(z))\ad y)}=0.
\end{equation}

\bigskip

The invariance under the adjoint action may be thought of as how the
Killing form behaves under an ``infinitesimal'' transformation; however,
recall that $\ad$ also induces the ``finite'' transformation (\ref{coact}).
Using this fact, one sees that the left-hand side of the second of
(\ref{inv-inf}) becomes
\begin{eqnarray}
\killing{\rho}{(z_{(1)}\ad x)}{(z_{(2)}\ad y)}&=&
\killing{\rho}{x^{(1)}}{y^{(1)}}\inprod{z_{(1)}}{x^{(2)'}}
\inprod{z_{(2)}}{y^{(2)'}}\nonumber\\
&=&\killing{\rho}{x^{(1)}}{y^{(1)}}\inprod{\Delta(z)}{x^{(2)'}
\otimes y^{(2)'}}\nonumber\\
&=&\killing{\rho}{x^{(1)}}{y^{(1)}}\inprod{z}{x^{(2)'}
y^{(2)'}}.
\end{eqnarray}
The right-hand side is equivalent to $\killing{\rho}{x}{y}\inprod{z}{\IA}$,
so since these must agree for all $z\in\U$,
\begin{equation}
\killing{\rho}{x^{(1)}}{y^{(1)}}x^{(2)'}y^{(2)'}=
\killing{\rho}{x}{y}\IA.\label{inv-fin}
\end{equation}
which expresses the ``finite'' invariance of the Killing form.

\subsection{The Killing Metric for a Quantum Lie Algebra}

In the case when $\U$ is not only a QHA but also a QLA with generators
$\{\gen{A}|A=1,\ldots,n\}$, one can define the {\em Killing metric}
associated with the rep $\rho$ as
\begin{equation}
\metric{\rho}{AB}:=\killing{\rho}{\gen{A}}{\gen{B}}.
\end{equation}

It is now convenient to introduce the quantities $\I{\rho}{A}:=
\trace{\rho}{\gen{A}}$; one can deduce from the results of the previous
sections that these ``deformed traces'' must satisfy
\begin{eqnarray}
\trace{\rho}{(\gen{A}\ad\gen{B})}&=&\struc{AB}{C}\I{\rho}{C}=0,
\nonumber\\
\trace{\rho}{(\O{A}{B}\ad\gen{C})}&=&\bigR{DB}{AC}\I{\rho}{D}=\delta^B_A
\I{\rho}{C},\nonumber\\
\trace{\rho}{\gen{B}}\bigA{B}{A}&=&\I{\rho}{B}\bigA{B}{A}=\I{\rho}{A}
\IA.\label{I-relations}
\end{eqnarray}
The first of (\ref{I-relations}) implies that by multiplying
(\ref{Lie-comm}) by $u$ and tracing over $\rho$, $\metric{\rho}{AB}$
satisfies
\begin{equation}
\metric{\rho}{AB}=\bigR{CD}{AB}\metric{\rho}{CD}=\bigD{C}{A}
\metric{\rho}{BC}
\end{equation}
(see (\ref{bigD-def}) for the definition of $\bigD{}{}$), which expresses
the ``symmetry'' of the Killing metric.  One can also obtain the ``total
antisymmetry'' of the structure constants in a similar way; since the
counits of all the generators vanish, (\ref{inv-inf}) requires that
\begin{eqnarray}
0&=&\killing{\rho}{(\gen{C(1)}\ad\gen{A})}{(\gen{C(2)}\ad\gen{B})}
\nonumber\\
&=&\killing{\rho}{(\gen{C}\ad\gen{A})}{\gen{B}}+\killing{\rho}{(\O{C}{D}
\ad\gen{A})}{(\gen{D}\ad\gen{B})}\nonumber\\
&=&\struc{CA}{D}\killing{\rho}{\gen{D}}{\gen{B}}+\bigR{ED}{CA}\struc{DB}{F}
\killing{\rho}{\gen{E}}{\gen{F}}.
\end{eqnarray}
Therefore,
\begin{equation}
\struc{CA}{D}\metric{\rho}{DB}+\bigR{ED}{CA}\struc{DB}{F}\metric{\rho}{EF}
=0.\label{asymm}
\end{equation}

By using (\ref{gens-coact}), together with (\ref{inv-fin}), the invariance
of the Killing metric under finite transformations takes the form
\begin{equation}
\metric{\rho}{CD}\bigA{C}{A}\bigA{D}{B}=\metric{\rho}{AB}\IA .
\label{inv-finite}
\end{equation}

\subsection{Quadratic Casimirs}

Now, suppose that $\metric{\rho}{AB}$ is invertible, \ie there exists a
numerical matrix $\invmet{\rho}{AB}$ such that
\begin{equation}
\metric{\rho}{AC}\invmet{\rho}{CB}=\invmet{\rho}{BC}\metric{\rho}{CA}=
\delta^B_A .
\end{equation}
Then (\ref{inv-finite}) implies that $\bigA{A}{C}\bigA{B}{D}
\invmet{\rho}{CD}=\invmet{\rho}{AB}\IA$; this in turn indicates that
the {\em quadratic casimir} defined by
\begin{equation}
\cas{\rho}:=\invmet{\rho}{AB}\gen{A}\gen{B}
\end{equation}
is central.  Why?  First, note that $\cas{\rho}$ is right-invariant:
\begin{eqnarray}
\DA (\cas{\rho})&=&\invmet{\rho}{AB}\DA (\gen{A})\DA (\gen{B})\nonumber\\
&=&\gen{C}\gen{D}\otimes\invmet{\rho}{AB}\bigA{C}{A}\bigA{D}{B}\nonumber\\
&=&\gen{C}\gen{D}\invmet{\rho}{CD}\otimes\IA\nonumber\\
&=&\cas{\rho}\otimes\IA .
\end{eqnarray}
Now, recall (\ref{coact}) and (\ref{adj-mult}); the first of these states
that if $y\in\U$ is right-invariant, $x\ad y=\epsilon (x)y$ for all $x$.
The second therefore implies
\begin{equation}
xy=(x_{(1)}\ad y)x_{(2)}=\epsilon(x_{(1)})yx_{(2)}=yx,
\end{equation}
namely, {\em any} right-invariant element of $\U$ is central.  Since the
right-invariance of $\cas{\rho}$ has just been demonstrated, it follows
that the quadratic casimir commutes with everything in the algebra, just as
in the classical case.

\subsection{The Canonical Killing Metric}\label{canonical}

Suppose that there exists at least one rep $\rho_0$ for which the
associated Killing metric $\metric{0}{AB}$ is invertible, \ie
$\invmet{0}{AB}$ exists.  Given another rep $\rho$, define the numerical
matrix $\K{A}{B}$ by
\begin{equation}
\K{A}{B}:=\invmet{0}{AC}\metric{\rho}{CB},\label{def-K}
\end{equation}
so that $\metric{\rho}{AB}\equiv\metric{0}{AC}\K{C}{B}$.  Therefore,
(\ref{asymm}) may be rewritten in terms of $\metric{0}{}$ and $\K{}{}$:
\begin{eqnarray}
0&=&\struc{CA}{D}\left( \metric{0}{DF}\K{F}{B}\right) +\bigR{ED}{CA}
\struc{DB}{F}\left( \metric{0}{EM}\K{M}{F}\right) \nonumber\\
&=&-\left( \bigR{ED}{CA}\struc{DF}{M}\metric{0}{EM}\right) \K{F}{B}+
\bigR{ED}{CA}\struc{DB}{F}\metric{0}{EM}\K{M}{F}\nonumber\\
&=&\bigR{ED}{CA}\metric{0}{EM}\left( \K{M}{F}\struc{DB}{F}-\struc{DF}{M}
\K{F}{B}\right) .
\end{eqnarray}
Since both $\bigR{}{}$ and $\metric{0}{}$ are invertible, this may be
multiplied by $\matrix{(\widehat{\real}^{-1})}{CA}{PN}\invmet{0}{LP}$ to
obtain the matrix equation
\begin{equation}
\matrix{\comm{\K{}{}}{\adj(\gen{N})}}{L}{B}=0,\label{K-comm}
\end{equation}
\ie the matrix $\K{}{}$ commutes with all the generators in the adjoint rep.
Therefore, if the adjoint rep is irreducible, Schur's lemma can be applied,
which would indicate that $\K{}{}$ is proportional to the identity matrix.
This in turn gives the result that $\metric{\rho}{AB}\propto\metric{0}{AB}$
for {\em all} $\rho$, so that, up to a multiplicative constant, all Killing
metrics may be written in terms of a canonical, \ie rep-independent,
Killing metric $\eta_{AB}$.  Once the normalization of $\eta$ is fixed,
this allows the definition of the {\em index} $\index{\rho}$ of the rep
$\rho$ in exact analogy with the classical case, \ie $\metric{\rho}{AB}
\equiv\index{\rho}\eta_{AB}$.

The existence of at least one rep where the Killing metric is invertible
has been assumed in the above proof, so since all metrics are proportional
to $\eta_{AB}$, the canonical Killing metric itself must be invertible.
This fact implies the existence of a rep-independent quadratic casimir $Q:=
\eta^{AB}\gen{A}\gen{B}$, which of course will still be central in the QLA.
Therefore, as in the classical case, for any rep $\rho$, $Q\equiv
\index{\rho}\cas{\rho}$.

\subsection{Comments on the Adjoint Representation}

Before proceeding to the particular case where the QLA of interest is
$\uea{su}{N}$, it is necessary to consider some of the consequences of the
trace relations (\ref{I-relations}), specifically as pertains to the
adjoint rep.  These relations are obviously trivially satisfied for
$\I{\rho}{A}\equiv 0$, which is the case for the classical compact Lie
algebras, where the generators are traceless in all reps.  However, there
is no reason to assume that this tracelessness condition still holds when
the object in question is a QLA (or even a classical {\em noncompact} Lie
algebra, for that matter), \ie there may exist a rep $\rho$ such that
$\I{\rho}{A}$ does not vanish.  If this is in fact the case, one is then
able to deduce the existence of another numerical object $\D{A}$ which
satisfies
\begin{equation}
\struc{AB}{C}\D{B}=0.\label{null}
\end{equation}
Why should this quantity exist?  From the last of (\ref{I-relations}),
$\I{\rho}{A}\IA$ is an algebra-valued eigenvector of $\bigA{t}{}$ with
eigenvalue unity.  The transpose of any matrix has the same eigenvalues as
the original, so this implies the existence of a numerical quantity $\D{A}$
such that $\D{A}\IA$ is a nonzero algebra-valued eigenvector of $\bigA{}{}$
with unit eigenvalue, \ie
\begin{equation}
\bigA{A}{B}\D{B}=\D{A}\IA.\label{D-defn}
\end{equation}
This in turn implies (\ref{null}), as well as
\begin{equation}
\bigR{CA}{BD}\D{D}=\delta^A_B\D{C}.
\end{equation}

(\ref{gens-coact}) and (\ref{D-defn}) together give the important result
that the quantity
\begin{equation}
\gen{0}:=\D{A}\gen{A}
\end{equation}
is right-invariant, and therefore commutes with the entire algebra.  Now,
assume that the deformed trace of this element
$\I{\rho}{0}:=\trace{\rho}{\gen{0}}$ vanishes for all reps.  Therefore, if
$\rho$ is an irrep of $\U$, then Schur's lemma states that there exists a
$\k$-number $\mu(\rho)$ such that $\rho(\gen{0})=\mu(\rho)I$, where $I$ is
the identity matrix in the corresponding irrep.  Thus, $\I{\rho}{0}=\mu(
\rho)\trace{\rho}{}$.  The quantity $\trace{\rho}{}$ may be thought of as
the ``deformed dimension'' of the irrep, and it is assumed here that it is
always nonvanishing; hence, the condition $\I{\rho}{0}=0$ implies that
$\mu(\rho)$ vanishes.  However, if this holds for {\em all} irreps, then
$\gen{0}$ must be identically zero as an element of $\U$.  The linear
independence of the $n$ generators $\{\gen{A}\}$ then leads to the
conclusion that $\D{A}\equiv 0$, contrary to the initial assumption, namely
that there exists a nonvanishing $\D{A}$.  Therefore, there must exist a
rep $\rho'$ for which $\pI{0}:=\I{\rho'}{0}\neq 0$, and hence the
individual traces $\pI{A}:=\I{\rho'}{A}$ are also nonvanishing.

The existence of this particular rep allows the definition of a new set of
generators $\{\pgen{A}|A=1,\ldots,n\}$ as
\begin{equation}
\pgen{A}:=\gen{A}-\frac{\pI{A}}{\pI{0}}\gen{0}.
\end{equation}
(Note that these are traceless in the rep $\rho'$.)  The
commutation relations and adjoint actions now take the forms
\begin{equation}
\pgen{A}\pgen{B}-\bigR{CD}{AB}\pgen{C}\pgen{D}=\left[ \struc{AB}{C}-
(\delta^D_A\delta^C_B-\bigR{CD}{AB})\frac{\pI{D}}{\pI{0}}
\gen{0}\right] \pgen{C}\label{comm-prime}
\end{equation}
and
\begin{eqnarray}
\gen{0}\ad\gen{0}=0,&&\pgen{A}\ad\gen{0}=0,\nonumber\\
\gen{0}\ad\pgen{A}=\D{B}\struc{BA}{C}\pgen{C},&&\pgen{A}\ad\pgen{B}=\left(
\delta_A^D-\frac{\pI{A}}{\pI{0}}\D{D}\right) \struc{DB}{C}\pgen{C},
\label{adj-prime}
\end{eqnarray}
However, the primed generators are {\em not} linearly independent, because
$\D{A}\pgen{A}\equiv 0$, so they do not form a basis for $\g$.  However, if
one of them, for instance $\pgen{n}$, is replaced by $\gen{0}$, then $\{
\gen{0},\pgen{a}|a=1,\ldots,n-1\}$ is a linearly independent collection of
$n$ vectors, and therefore is a basis for $\g$.  (For the rest of this
work, when a capital index $A$ appears on a primed generator, it will take
the values $A=0,1,\ldots,n-1$, with $\pgen{0}=\gen{0}$, and a small index
$a$ the values $a=1,\ldots,n- 1$.)

The structure constants $f'$ in this new basis, and thus the
matrices in the adjoint rep, are easily obtained from (\ref{adj-prime}),
and the reader can immediately see that
\begin{equation}
\pstruc{00}{0}=\pstruc{00}{a}=\pstruc{a0}{0}=\pstruc{a0}{b}=\pstruc{0a}{0}
=\pstruc{ab}{0}=0,\label{struc-prime}
\end{equation}
\ie the only nonvanishing structure constants are $\pstruc{Aa}{b}$.  In
this basis all the matrices in the adjoint rep will be block-diagonal, with
a zero in the upper left-hand corner and an $(n-1)\times (n-1)$ matrix as
the lower right-hand submatrix (corresponding to the $(00)$ and $(ab)$
entries respectively).  This set of matrices closes under multiplication,
so the $(n-1)\times (n-1)$ matrices
\begin{equation}
\matrix{\padj{}}{a}{b}(\pgen{A}):=\pstruc{Ab}{a}
\end{equation}
form a rep of the QLA.

\section{The Quantum Lie Algebra $\uea{su}{N}$}
\setcounter{equation}{0}

The discussions and conclusions just presented in
Section \ref{chap-Killing} closely parallel those usually found when
talking about undeformed compact Lie algebras, \eg the existence of a rep
for which the associated Killing metric is invertible, irreducibility of
the adjoint rep, tracelessness of the basis, \etc For the classical compact
Lie algebras, all these assumptions hold: the Killing metric in the adjoint
rep is positive definite (this is actually the {\em definition} of a
compact Lie algebra) and thus invertible, the adjoint rep is an irrep, and
all the generators are traceless (the ``$S$'' in $SU(N)$).  Thus, results
like the existence of $\eta_{AB}$ and a rep-independent quadratic casimir
follow.  However, these assumptions may not be valid for a {\em deformed}
Lie algebra; hence, to use any of the conclusions of the previous section,
one must first ask how many of these still hold for a QLA.

Because the classical Lie algebra {\goth su}$(N)$ is the one which is most
familiar to most physicists, the deformed version $\uea{su}{N}$ is the best
example with which to illustrate many of the results just obtained.  This
Section will examine this particular QLA in detail.

\subsection{R-Matrix Construction}\label{funad}

As shown in Section \ref{QQLA}, if a numerical R-matrix for a QHA is given,
a QLA may be constructed; $\uea{su}{N}$ is the QLA found by this method
using the numerical R-matrix for the QG $SU_q(N)$, which in the fundamental
rep is given by multiplying the R-matrix for $A_{N-1}$ in \cite{RTF} by
$q^{-\inv{N}}$ \cite{linear}:
\begin{equation}
R=q^{-\inv{N}}\left( q\sum_IE_{II}\otimes E_{II}+\sum_{I\neq J}E_{II}
\otimes E_{JJ}+\lambda\sum_{I>J}E_{IJ}\otimes E_{JI}\right) ,
\end{equation}
where $E_{IJ}$ is the $N\times N$ numerical matrix whose only nonzero entry
is a 1 at $(I,J)$, the tensor product which appears is that between
numerical spaces, and $\k=\real$.  ($\lambda:=q-q^{-1}$ is the parameter
referred to earlier which describes how the QG differs from the classical
group as a function of $q$.)  The normalization for the matrix $D$ is
chosen so that $D=\mbox{diag}(1,q^2,\ldots,q^{2(N-1)})$, which fixes the
constants $\alpha$ and $\beta$ from the Appendix to be $q^{2N-1-\inv{N}}$
and $q^{1-\inv{N}}$ respectively.

The commutation relations (\ref{X-X}) and adjoint actions (\ref{X-ad-X})
follow, with the structure constants given by (\ref{struc-R}).  Now that
they are known explicitly in terms of the above numerical R-matrix, one can
immediately obtain several results specific to $\uea{su}{N}$: first of all,
suppose that $V_{(ij)}:=\matrix{V}{i}{j}$ is a null eigenvector of the
transposed adjoint matrices, \ie
\begin{equation}
\struc{(ij)(k\ell )}{(rs)}V_{(rs)}\equiv\inv{\lambda}
\matrix{\widetilde{R}}{mk}{jn}\matrix{(\Rhat V_2-\Rhat^{-1}V_2
\Rhat^2)}{ni}{m\ell}=0.
\end{equation}
{}From this, it follows that $V_2$ must commute with $\Rhat^2$.  This implies
that $V$ must be a multiple of the identity matrix: to see why, suppose
$V_2$ commutes with a matrix $\widehat{M}$ such that $\widetilde{M}$
exists, where $\matrix{M}{ij}{k\ell}:=\matrix{\widehat{M}}{ji}{k\ell}$.  If
this is the case, then
\begin{eqnarray}
0&=&\matrix{\widetilde{M}}{mj}{kn}\matrix{(\widehat{M}V_2-V_2\widehat{M})
}{ni}{m\ell}\nonumber\\
&=&\matrix{\widetilde{M}}{mj}{kn}\matrix{M}{in}{mr}\matrix{V}{r}{\ell}-
\matrix{\widetilde{M}}{mj}{kn}\matrix{V}{i}{s}\matrix{M}{sn}{m\ell}
\nonumber\\
&=&\delta^i_k\matrix{V}{j}{\ell}-\matrix{V}{i}{k}\delta^j_{\ell},
\end{eqnarray}
which can only be satisfied if $\matrix{V}{i}{j}\propto\delta^i_j$.  The
numerical R-matrix for $SU_q(N)$ satisfies the quadratic characteristic
equation $\Rhat^2-q^{-\inv{N}}\lambda\Rhat-q^{-\frac{2}{N}}I=0$, so $V_2$
commutes with $\Rhat$.  Since $\widetilde{R}$ exists, the above result
applies, and all null eigenvectors of the transposed adjoint matrices must
be proportional to the identity.  This implies that since the traces
$\I{\rho}{(ij)}$ are such eigenvectors, all of them are equal to a
rep-dependent constant times the identity.  (From (\ref{QLA-I}), it is
easily seen that these also satisfy the first two of (\ref{I-relations}) as
well.)  In fact, for the fundamental rep $\fun$, in which the matrices take
the form
\begin{equation}
\matrix{\fun}{i}{j}(\gen{(k\ell)})=\inv{\lambda}\matrix{(I-\Rhat^2)
}{ik}{j\ell},\label{fn}
\end{equation}
a quick calculation gives
\begin{equation}
\I{\fun}{(ij)}=q^{-\inv{N}}\left( \quint{\inv{N}}{q}\quint{N}{q^{-1}}-1
\right) \delta^i_j
\end{equation}
(where I have used the standard notation for the ``quantum number'',
\begin{equation}
\quint{m}{q}:=\frac{q^{2m}-1}{q^2-1},
\end{equation}
so called because as $q\rightarrow 1$, $\quint{m}{q}\rightarrow m$).  In
the classical limit, this vanishes for all $N$.  This should come as no
surprise: after all, the adjoint rep of the classical Lie algebra {\goth
su}$(N)$ is irreducible, so there cannot exist any common nonzero null
eigenvector for the matrices $\adj(\gen{A})^t$.  However, for $q\neq 1$,
this is nonvanishing, and so it implies the existence of $\gen{0}$ and will
allow the construction of the $N^2-1$ primed generators, which together
with $\gen{0}$ give a basis for $\g$.  The same sort of argument allows one
to determine that all matrices $W^{(ij)}:=\matrix{W}{j}{i}$ satisfying
$\struc{(ij)(k\ell)}{(rs)}W^{(k\ell )}=0$ must be proportional to
$\matrix{(D^{-1})}{j}{i}$.  From (\ref{null}), this includes $\D{A}$, so
${\cal D}^{(ij)}\propto\matrix{(D^{-1})}{j}{i}$.

The arguments just given have two immediate consequences: first of all,
assume there exists a vector $\vec{e} =(e^1,\ldots,e^{n-1})$ which
satisfies $\padj(\pgen{A})\cdot\vec{e}=\vec{0}$; since each
$\adj(\pgen{A})$ is block-diagonal with a zero as its $(00)$ entry, then
any vector of the form $(\nu,\vec{e})$, $\nu$ being some constant, will be
a null eigenvector of $\adj(\pgen{A})$.  However, it has just been shown
that all such eigenvectors are proportional to $D^{-1}$, and thus to
$\D{A}$.  In the primed basis, $\D{0}=1$ and $\D{a}=0$, so any vector
annihilated by all the adjoint matrices must be a multiple of
$(1,\vec{0})$.  This implies that $\vec{e}$ must vanish, and that there are
{\em no} shared nonzero null eigenvectors of the matrices
$\padj(\pgen{A})$.  One therefore concludes that the $(N^2-1)$-dimensional
rep of $\uea{su}{N}$ given by $\padj$ is in fact an irrep, and
$\padj(\gen{0})$, being central, will be proportional to the identity
matrix.

Second, notice that
\begin{equation}
\frac{\I{\rho}{(ij)}}{\I{\rho}{0}}{\cal D}^{(k\ell )}=\inv{\tr(D^{-1})}
\delta^i_j\matrix{(D^{-1})}{\ell}{k},\label{rep-ind}
\end{equation}
which is {\em independent} of the constants of proportionality for both
$\I{\rho}{(ij)}$ and ${\cal D}^{(ij)}$.  The generators $\gen{0}$ and
$\pgen{(ij)}$ therefore take the forms
\begin{eqnarray}
\gen{0}=\tr(D^{-1}X),&&\pgen{(ij)}=\matrix{X}{i}{j}-\inv{\tr(D^{-1})}
\tr(D^{-1}X)\delta^i_j.\label{def-su}
\end{eqnarray}
(The above expression for $\gen{0}$ together with (\ref{inv-tr}) provides
further proof that $\gen{0}$ is right-invariant, and thus central.)  A
consequence of the rep-independence of (\ref{rep-ind}) is that the primed
generators are traceless in {\em all} reps, not just the one for which
$\pI{0}$ is nonvanishing (\eg $\fun$).  Thus, when one of them is replaced
by $\gen{0}$ to obtain a linearly independent basis for $\g$, the result is
that the QLA is generated by the $(N^2-1)$-dimensional subspace of
traceless basis elements $\g'=\{\pgen{a}\}$ and a central element
$\gen{0}$.  In fact, this also implies that $\gen{0}$ must vanish in the
classical limit, because $\g'=\mbox{\goth su}(N)$ in this case, and one
requires that $\uea{su}{N}\rightarrow U(\mbox{\goth su}(N))$.

\bigskip

$\uea{su}{N}$ supposedly describes a deformation of a unitary algebra, so a
hermiticity condition on its elements must be imposed: this is given in
terms of the involution acting on the matrices $L^{\pm}$ as $\theta(L^{\pm}
):=S(L^{\mp})^t$, which implies that the generator matrix is hermitian, \ie
$\theta(X)=X^t$, or $\theta(\gen{(ij)}) =\gen{(ji)}$.  If the subspace $\g$
is to be a vector space over $\real$ as in the classical case, then for
$x\in\g$, $\theta(x)=x$.  Being just a linear combination of the
generators, such an element $x$ may therefore be written
\begin{equation}
x=\xi^A\gen{A}=\tr(\Xi X),\label{x-Xi}
\end{equation}
where $\xi^{(ij)}\equiv\matrix{\Xi}{j}{i}$ are the numerical coordinates of
$x$ in the basis $\{\gen{(ij)}\}$.  For $x$ to be self-adjoint, the matrix
$\Xi$ must be hermitian.

\subsection{Properties of the Killing Metric}\label{killsu}

In the fundamental rep $\fun$, the generators are given by (\ref{fn}), and
$u$ by $q^{1+\inv{N}-2N}D$ (see Appendix), so the Killing metric associated
with the fundamental rep in the basis $\{\gen{(ij)}\}$ may be found
explicitly:
\begin{eqnarray}
\metric{\fun}{(ij)(k\ell)}&=&q^{-\inv{N}}\left( q\quint{1-\inv{N}}{q}-
\inv{q}\quint{1+\inv{N}}{q^{-1}}+q^{\frac{2}{N}-3}\quint{\inv{N}}{
q^{-1}}^2\quint{N}{q^{-1}}\right) \delta^i_j\delta^k_{\ell}\nonumber\\
&&+q^{1-\frac{3}{N}-2N}\delta^i_{\ell}\matrix{D}{k}{j}.\label{kill-sun}
\end{eqnarray}

Now, suppose one switches to the primed basis and computes the Killing
metric for a rep $\rho$ which reduces to irreps $\{\rho_i\}$.  Then
\begin{eqnarray}
\metric{\rho}{0a}&=&\trace{\rho}{\gen{0}\pgen{a}}\nonumber\\
&=&\sum_i\trace{\rho_i}{\gen{0}\pgen{a}}\nonumber\\
&=&\sum_i\mu(\rho_i)\trace{\rho_i}{\pgen{a}},
\end{eqnarray}
which vanishes due to the tracelessness of the primed generators in all
reps.  Furthermore, since $\gen{0}$ is central, then $\metric{\rho}{a0}
\equiv\metric{\rho}{0a}$, so the Killing metric will be a block-diagonal
matrix with $\metric{\rho}{00}$ in the upper left-hand corner and the
$(N^2-1)\times (N^2-1)$ matrix $\metric{\rho}{ab}$ in the lower right-hand
corner.

Any element $x\in\g$ may be written either as in (\ref{x-Xi}), or in terms
of the primed basis as $x=\xi^0\gen{0}+\xi^a\pgen{a}$, from which it
follows from (\ref{def-su}) that $\xi^0=\frac{\tr\Xi}{\tr D^{-1}}$.  As the
reader has just seen, the Killing metric is block-diagonal in the primed
basis, so $\metric{\fun}{AB}\xi^A\xi^B=\metric{\fun}{00}(\xi^0)^2+
\metric{\fun}{ab}\xi^a\xi^b$.  Therefore, by using (\ref{kill-sun}),
\begin{eqnarray}
\metric{\fun}{ab}\xi^a\xi^b&=&\metric{\fun}{AB}\xi^A\xi^B-\metric{\fun}{00}
(\xi^0)^2\nonumber\\
&=&\metric{\fun}{(ij)(k\ell)}\left[ \matrix{\Xi}{j}{i}\matrix{\Xi}{\ell}{k}
-\matrix{(D^{-1})}{j}{i}\matrix{(D^{-1})}{\ell}{k}\left(\frac{\tr\Xi}{\tr
D^{-1}}\right) ^2\right]\nonumber\\
&=&q^{1-\frac{3}{N}-2N}\left[ \tr (D\Xi^2)-\frac{(\tr\Xi)^2}{\tr
D^{-1}}\right]\nonumber\\
&=&q^{1-\frac{3}{N}-2N}\tr\left[ D\left( \Xi-\frac{\tr\Xi}{\tr D^{-1}}
D^{-1}\right) ^2\right]\nonumber\\
&=&q^{1-\frac{3}{N}-2N}\tr\abs{D^{\half}\Xi-\frac{\tr\Xi}{\tr D^{-1}}D^{
-\half}}^2,
\end{eqnarray}
where the fact that $D$ is hermitian and positive definite (recall that
there is a basis where it is diagonal with each entry being the square of a
real number) has been used.  The quantity above is thus always nonnegative,
and vanishes only when $\Xi\propto D^{-1}$, \ie $x\propto\gen{0}$ and
$\xi^a=0$.  The Killing metric $\metric{\fun}{ab}$ constructed using only
the primed generators is positive definite, and thus invertible.

Therefore, the matrix $\K{}{}$ from (\ref{def-K}) can be constructed by
using $\rho_0=\fun$, and it too will be block-diagonal.  In particular, its
upper left-hand entry will be $\invmet{\fun}{00}\metric{\rho}{00}$.
However, when $\K{}{}$ is then multiplied by any matrix in the adjoint rep,
the result will have a zero in this spot, since all adjoint matrices have a
vanishing entry there.  The commutation relation (\ref{K-comm}) hence says
nothing about this particular entry, and instead becomes a statement only
about the $(N^2-1)\times (N^2-1)$ submatrices.  In fact, since it has
already been demonstrated that the rep of the algebra given by the
$(N^2-1)\times (N^2-1)$ submatrices of the adjoint is an irrep, the result
obtained in Section \ref{canonical} follows, and all of the submetrics
$\metric{\rho}{ab}$ will be proportional to some canonical Killing metric
$\eta_{ab}$.  Since this includes the invertible metric
$\metric{\fun}{ab}$, it follows that the canonical Killing metric is itself
invertible as well.

This fact implies that the rep-independent quadratic casimir $Q$ can be
constructed, and therefore, together with $\gen{0}$, classifies irreps of
$\uea{su}{N}$.  However, due to the block-diagonality of the Killing
metric, it is immediately seen that the rep-independent quadratic casimir
computed using only the primed generators, \ie $\pcas:=\eta^{ab}\pgen{a}
\pgen{b}$, may be written as
\begin{equation}
\pcas=Q-\eta^{00}\gen{0}{}^2.
\end{equation}
Both terms on the right-hand side are central, so $\pcas$ is as well.
Therefore, in classifying an irrep $\rho$ of $\uea{su}{N}$, $Q$ may be
replaced by $\pcas$.

\subsection{Hopf Algebraic Structure}

The primed generators are just linear combinations of the original ones, so
they all have vanishing counit.  Also, (\ref{struc-prime}) indicates that
$\pstruc{Aa}{0}=0$, and therefore $\U\ad\g'\subseteq\g'$.  $\g'$ hence is a
subspace of $\U$ which satisfies criteria (II) and (IV) of the definition
of a QLA.  Furthermore, as was shown in \ref{funad}, it is spanned by
traceless generators and its $(N^2-1)$-dimensional adjoint rep is an irrep.
These facts lead to the results in Section \ref{killsu}, in which it was
found that the Killing metric and quadratic casimir computed using only the
elements of $\g'$ have all the properties which one would expect for the
analagous quantities for the undeformed Lie algebra {\goth su}$(N)$, \ie
invertibility, centrality, \etc Also, $\gen{0}\rightarrow 0$ in the
classical limit $\uea{su}{N}\rightarrow U(\mbox{\goth su}(N))$.

All of this certainly seems to suggest that it is possible to follow the
classical example, and define the QLA $\uea{su}{N}$ to be that subspace of
$\U$ generated by $\g'$.  Unfortunately, this isn't the case; if the reader
recalls (\ref{comm-prime}), s/he sees that the primed generators do {\em
not} close under multiplication--$\gen{0}$ makes an appearance.  It can
also be shown that $\gen{0}$ appears explicitly in the coproducts of the
primed generators as well, so $U_q(\g')$ doesn't satisfy (I) or (III).
$\gen{0}$ must be included as a generator, and therefore one cannot follow
the classical example and obtain {\goth su}$(N)$ from {\goth gl}$(N)$
simply by throwing away the ``traceful'' central element after imposing the
hermiticity condition.  (In contrast, at the QG level, there is a way to
obtain $SU_q(N)$ and preserve the HA structure \cite{linear}.)

However, from the point of view of an irrep $\rho$ of $\uea{su}{N}$, the
exact HA structure is not so important; the coproduct, counit and antipode
will not enter into the discussion, and $\gen{0}$ may be replaced by
$\mu(\rho)$ in relations like (\ref{comm-prime}).  The other good quantum
numbers will be the ones obtained from considering $\g'$, namely $\pcas$
and the weights associated with the Cartan subalgebra of $\g'$, as in the
classical case.  Of course, if one wants to find higher reps from the
tensor product of two irreps, the HA structure is absolutely vital, and
will generally give results different from the classical case (\eg if
$N=2$, and $\rho_j$ is the spin $j$ irrep, $\rho_{j_1}\otimes\rho_{j_2}$
may no longer be $\rho_{\abs{j_1-j_2}}\oplus\ldots\oplus\rho_{j_1+j_2}$,
but some deformed version instead).

\subsection{Other Quantum Lie Algebras}

Recall a very important assumption made in the discussion of $\uea{su}{N}$,
namely, that the numerical R-matrix satisfies a quadratic characteristic
equation.  This lead to the result that all of the traces $\I{\rho}{A}$
were proportional to one another, and thus the primed generators (except
$\gen{0}$, of course) were traceless in all reps.  It also lead to an
explicit form of $\D{A}$, which then implied the irreducibility of the rep
$\padj$.  Once these were in hand, the canonical Killing metric and the
like followed.

For the QG $SU_q(N)$, this was in fact the case; however, what about the
other QLAs corresponding to the deformed versions of the classical compact
Lie algebras {\goth so}$(N)$ and {\goth sp}$(\frac{N}{2})$?  As given in
\cite{RTF}, the numerical R-matrices of $SO_q(N)$ and $SP_q(\frac{N}{2})$
satisfy the {\em cubic} characteristic equations
\begin{equation}
(\Rhat-qI)(\Rhat+q^{-1}I)(\Rhat-\epsilon q^{\epsilon-N})=0,
\end{equation}
where $\epsilon=1$ for $SO_q(N)$ and $\epsilon=-1$ for $SP_q(\frac{N}{2})$.
Therefore, the condition that $V_2$ and $\Rhat^2$ commute does {\em not}
imply that $V_2$ commutes with $\Rhat$.

However, recall that the existence of a quadratic characteristic equation
was not essential; what was necessary was the existence of $\widetilde{M}$,
which for an arbitrary QLA means that one must be able to find the matrix
$\widetilde{(RPR)}$.  At the present time, the author does not know if this
matrix always exists for QLAs other that $\uea{su}{N}$; however, by using
(\ref{fn}) (which still gives the fundamental rep if the appropriate
numerical R-matrix is plugged in), $\I{\fun}{A}$ may be calculated for
these other QLAs, with the results being
\begin{equation}
\I{\fun}{(ij)}=\left( q^{\epsilon-N}-q^{N-\epsilon}\right)\delta^i_j.
\end{equation}
So, even though it has not been shown for {\em all} reps, this at least
seems to hint that the same results may be obtained for $\uea{so}{N}$ and
$\uea{sp}{\frac{N}{2}}$ as for $\uea{su}{N}$, and therefore the same
conclusions follow.

This is not the whole story, however; for $\uea{su}{N}$, an irreducible
$(N^2-1)$-dimensional adjoint rep was easily interpreted as a deformation
of the classical case.  However, the adjoint reps for {\goth so}$(N)$ and
{\goth sp}$(\frac{N}{2})$ have dimensions given by $\half N(N-\epsilon)$,
due to the fact that both are endowed with a metric $C$ and an accompanying
orthogonality restriction on their elements.  This metric also exists in
the deformed case, and obviously must be taken into account when discussing
the adjoint rep.  The QLA of these groups has been studied (\cf
\cite{WSWW}), but how to extend the results of this work to such QLAs
remains an open problem.

\section{Example:  $\uea{su}{2}$}
\setcounter{equation}{0}

In this Section, some explicit calculations for the case of the QLA
$\uea{su}{2}$ will be presented, and will (hopefully) serve to illustrate
many of the results from the previous sections.

Consider the unital associative algebra generated by the elements $\{H,
X_+,X_-\}$, modulo the Jimbo-Drinfel'd commutation relations
\cite{Jimbo,Drin}
\begin{eqnarray}
\comm{H}{X_{\pm}}&=&\pm 2X_{\pm},\nonumber\\
\comm{X_+}{X_-}&=&\frac{q^H-q^{-H}}{\lambda},
\end{eqnarray}
where $q\in\real$.  This algebra is actually a QHA: the coproducts,
counits, antipodes, and conjugates of these elements are given by
\begin{eqnarray}
\Delta(H)=H\otimes\IU+\IU\otimes H,&&\Delta(X_{\pm})=X_{\pm}\otimes
q^{\half H}+ q^{-\half H}\otimes X_{\pm},\nonumber\\
\epsilon(H)=\epsilon(X_{\pm})=0,&&\nonumber\\
S(H)=-H,&&S(X_{\pm})=-q^{\pm 1}X_{\pm},\nonumber\\
\theta(H)=H,&&\theta(X_{\pm})=X_{\mp},
\end{eqnarray}
and the universal R-matrix by \cite{Rosso}
\begin{equation}
\R=\sum_{n=0}^{\infty}\frac{(1-q^{-2})^n}{\quint{n}{q}!}q^{\half(H\otimes
H+nH\otimes\IU-n\IU\otimes H)}X_+^n\otimes X_-^n ,
\end{equation}
where the ``quantum factorial'' is defined as
\begin{equation}
\quint{n}{q}!:=\left\{ \begin{array}{ll}
1&n=0,\\
\prod_{m=1}^n\quint{m}{q}&n=1,2,\ldots
\end{array}\right.
\end{equation}
This QHA is usually denoted by $\uea{su}{2}$, and is the deformed version
of the classical {\goth su}$(2)$ Lie algebra.

The fundamental reps for both the deformed and undeformed cases coincide,
\ie the matrices
\begin{eqnarray}
\fun(H)=\left( \begin{array}{cc}
-1&0\\0&1\end{array}\right) ,&
\fun(X_+)=\left( \begin{array}{cc}
0&0\\-1&0\end{array}\right) ,&
\fun(X_-)=\left( \begin{array}{cc}
0&-1\\0&0\end{array}\right) ,\label{fundrep}
\end{eqnarray}
satisfy the Jimbo-Drinfel'd commutation relations for any value of $q$.
For this rep, one can obtain the numerical R-matrix
\begin{equation}
R=q^{-\half}\left( \begin{array}{cccc}
q&0&0&0\\0&1&0&0\\0&\lambda &1&0\\0&0&0&q
\end{array}\right) ,\label{R-sl2}
\end{equation}
and the matrices $\L^{\pm}$ from (\ref{L-defs}):
\begin{eqnarray}
L^+=\left( \begin{array}{cc}
q^{-\half H}&-q^{-\half}\lambda X_+\\0&q^{\half H}
\end{array}\right) ,&&
L^-=\left( \begin{array}{cc}
q^{\half H}&0\\q^{\half}\lambda X_-&q^{-\half H}
\end{array}\right) .
\end{eqnarray}

One may now find the QLA from this QHA by taking the above expressions for
$L^{\pm}$ and constructing the $2\times 2$ matrix $X$.  The generators
$\gen{1}$, $\gen{+}$, $\gen{-}$ and $\gen{2}$ are defined to be the entries
of $X$, via
\begin{equation}
X=\left( \begin{array}{cc}
\gen{1}&\gen{+}\\
\gen{-}&\gen{2}
\end{array} \right) .\label{chi-fund}
\end{equation}
(Notice from the hermiticity condition that $\gen{1}$ and $\gen{2}$ are
self-adjoint, and $\theta(\gen{\pm})=\gen{\mp}$.)

Using the change of basis introduced earlier, it is easily seen that
$\gen{0}=\gen{1}+\inv{q^2}\gen{2}$ and $\pgen{\pm}=\gen{\pm}$.  However,
the connection with the classical case is a bit more obvious if, instead of
$\pgen{1}$, one defines $\gen{3}:=q^2\pgen{1}=\inv{\quint{2}{q^{-1}}}
(\gen{1}-\gen{2})$.  Using this basis, the adjoint actions are
\begin{eqnarray}
\gen{0}\ad\gen{0}=0,&\gen{a}\ad\gen{0}=0,&\gen{0}\ad\gen{a}=-\lambda
\quint{2}{q^{-1}}\gen{a}
\end{eqnarray}
(where $a=+,-,3$), as well as
\begin{eqnarray}
\gen{3}\ad\gen{3}=-\lambda\gen{3},&&\gen{\pm}\ad\gen{\mp}=\pm
\frac{\quint{2}{q^{-1}}}{q}\gen{3},\nonumber\\
\gen{3}\ad\gen{\pm}=\pm q^{\mp 1}\gen{\pm},&&\gen{\pm}\ad\gen{3}=\mp q^{\pm
1}\gen{\pm}.\label{adj-basis}
\end{eqnarray}

The commutation relations may also be found, albeit with a bit more work,
since $\bigR{}{}$ is a $16\times 16$ matrix.  However, when all is said and
done, one finds that $\gen{0}$ is indeed central as expected, and the other
generators satisfy
\begin{eqnarray}
q^{\mp 1}\gen{3}\gen{\pm}-q^{\pm 1}\gen{\pm}\gen{3}&=&\pm\left( \IU-\frac{
\lambda}{\quint{2}{q^{-1}}}\gen{0}\right) \gen{\pm},\nonumber\\
\gen{+}\gen{-}-\gen{-}\gen{+}&=&\frac{\quint{2}{q^{-1}}}{q}\left(
\IU-\frac{\lambda}{\quint{2}{q^{-1}}}\gen{0}\right) \gen{3}+\frac{\lambda
\quint{2}{q^{-1}}}{q}\gen{3}^2.\label{su2-comm}
\end{eqnarray}
(Note that if $\gen{0}$ is replaced by $-\lambda\quint{\half}{q}\quint{3
\over 2}{q^{-1}}I$, and then the redefinitions
\begin{eqnarray}
q\rightarrow q^{-\half},&\gen{\pm}\rightarrow\sqrt{q^{\half}\frac{
\quint{2}{q}}{\quint{2}{q^{\half}}}}\gen{\pm},&\gen{3}\rightarrow
\inv{\sqrt{\quint{2}{q}}}\gen{0},\label{subs}
\end{eqnarray}
are made, then one recovers the commutation relations for $SU_q(2)$ in
\cite{witten}.  As the reader will see in the next subsection, $\gen{0}$
does indeed have this value in the fundamental irrep.)

In the following two subsections, the two obvious reps to consider, namely,
the fundamental and adjoint, will be examined in detail.

\subsection{Fundamental Representation}

The numerical matrices for the generators $\{\gen{1},\gen{+},\gen{-},
\gen{2}\}$ in the $2\times 2$ fundamental rep $\fun$ of $\uea{su}{2}$ may be
found by using (\ref{fn}), where $R$ is given by (\ref{R-sl2}).  Then, by
switching the basis as described above, one obtains the matrices
\begin{eqnarray}
\fun(\gen{0})=-\lambda\quint{\half}{q}\quint{3\over 2}{q^{-1}}
\left( \begin{array}{cc}
1&0\\0&1\end{array}\right) ,&&
\fun(\gen{3})=\inv{\quint{2}{q^{-1}}}\left( \begin{array}{cc}
-1&0\\0&\inv{q^2}\end{array}\right) ,\nonumber\\
\fun(\gen{+})=\left( \begin{array}{cc}
0&0\\-\inv{q}&0\end{array}\right) ,&&\fun(\gen{-})=\left(
\begin{array}{cc}
0&-\inv{q}\\0&0\end{array}\right) .
\end{eqnarray}
(Note that as previously argued, $\gen{0}$ does in fact vanish as
$q\rightarrow 1$.)  It is of course necessary to find the matrix $D$, \ie
$u$ in this rep.  By using the explicit form of $R$ and the material in the
Appendix, this is a simple calculation, with the result being
\begin{equation}
\fun(u)=q^{-\frac{5}{2}}\left( \begin{array}{cc}
1&0\\0&q^2\end{array}\right) .
\end{equation}

The $4\times 4$ Killing metric using all the generators will be
block-diagonal, with the $(00)$ entry being $\metric{\fun}{00}=q^{-\half}
\quint{2}{q^{-1}}\quint{\half}{q}^2\quint{3\over 2}{q^{-1}}^2$ and the
$3\times 3$ metric computed in the basis $\{\gen{+},\gen{-},\gen{3}\}$
given by
\begin{equation}
\metric{\fun}{ab}=q^{-\frac{7}{2}}\left(
\begin{array}{ccc}
0&q&0\\\inv{q}&0&0\\0&0&\frac{q}{\quint{2}{q^{-1}}}\end{array}\right) .
\label{killing-fund}
\end{equation}
This must be proportional to the canonical Killing metric, of course, but
one is free to choose the normalization (and thus define the index).  The
choice in this work will be
\begin{equation}
\eta_{ab}:=\left( q+\inv{q}\right) \left( \begin{array}{ccc}
0&q&0\\\inv{q}&0&0\\0&0&\frac{q}{\quint{2}{q^{-1}}}\end{array}\right) ,
\end{equation}
so that $\index{\fun}=\frac{q^{-\frac{9}{2}}}{\quint{2}{q^{-1}}}$.  The
rep-independent quadratic casimir thus has the form
\begin{equation}
\pcas=\inv{q^2\quint{2}{q^{-1}}}\left( q^2\gen{+}\gen{-}+\gen{-}\gen{+}+
\quint{2}{q^{-1}}\gen{3}^2 \right) .
\end{equation}
(If the substitutions (\ref{subs}) are once again made, then this agrees
with the $SU_q(2)$ casimir in \cite{witten} up to an overall factor of
$\frac{q^2+ 1}{q(q+1)^2}$.)  When the matrix reps for the generators are
put into this expression, one finds
\begin{equation}
\fun(\pcas)=\frac{\quint{3}{q^{-1}}}{\quint{2}{q}\quint{2}{q^{-1}}}\left(
\begin{array}{cc}
1&0\\0&1\end{array}\right) .
\end{equation}

Recall that the classical expressions for the index and quadratic casimir
in the spin $j$ irrep are $\inv{3}j(j+1)(2j+1)$ and $j(j+1)I$ respectively,
and so for the fundamental rep, \ie $j=\half$, $\index{\fun}=\half$ and
$\fun(\pcas)=\frac{3}{4}I$.  For $q=1$, the reader can immediately see that
the results just obtained agree precisely with these values.

\subsection{Adjoint Representation}

Using the structure constants for $\uea{su}{2}$ from (\ref{adj-basis}), it
is trivial to find the generators in the $3\times 3$ adjoint rep $\padj$;
in the basis $\{\gen{+},\gen{-},\gen{3}\}$, they take the forms
\begin{eqnarray}
\padj(\gen{0})=-\lambda\quint{2}{q^{-1}}\left( \begin{array}{ccc}
1&0&0\\0&1&0\\0&0&1\end{array}\right) ,&&\padj(\gen{3})=\left(
\begin{array}{ccc}
\inv{q}&0&0\\0&-q&0\\0&0&-\lambda\end{array}\right) ,\nonumber\\
\padj(\gen{+})=\left( \begin{array}{ccc}
0&0&-q\quint{2}{q^{-1}}\\0&0&0\\0&\inv{q}&0\end{array}\right) ,&&
\padj(\gen{-})=\left( \begin{array}{ccc}
0&0&0\\0&0&\inv{q}\quint{2}{q^{-1}}\\-\inv{q}&0&0\end{array}\right) ,
\end{eqnarray}
and, again using the methods from the Appendix, $u$ is
\begin{equation}
\padj(u)=\left( \begin{array}{ccc}
\inv{q^2}&0&0\\0&\inv{q^6}&0\\0&0&\inv{q^4}\end{array}\right) .
\end{equation}
$\metric{\padj}{00}=\inv{q^2}\lambda^2\quint{2}{q^{-1}}^2\quint{3}{q^{-1}}$
in this case, and the $3\times 3$ Killing metric is
\begin{equation}
\metric{\padj}{ab}=\frac{\quint{4}{q^{-1}}}{q^3}\left( \begin{array}{ccc}
0&q&0\\\inv{q}&0&0\\0&0&\frac{q}{\quint{2}{q^{-1}}}\end{array}\right)
\label{killing-adj}
\end{equation}
so that $\index{\padj}=\frac{\quint{4}{q^{-1}}}{q^4\quint{2}{q^{-1}}}$.
The rep-independent quadratic casimir in the adjoint rep is
\begin{equation}
\padj(\pcas)=\frac{\quint{4}{q^{-1}}}{\quint{2}{q^{-1}}}\left(
\begin{array}{ccc}
1&0&0\\0&1&0\\0&0&1\end{array}\right) .
\end{equation}
Once again, for $q=1$, both $\index{\padj}$ and $\padj(\pcas)$ agree with
the classical values in the adjoint rep ($j=1$).

(The author suspects that the general forms for the two central quantities
$\gen{0}$ and $\pcas$ in the spin $j$ irrep $\rho_j$ are
\begin{eqnarray}
\rho_j(\gen{0})=-\lambda\quint{j}{q}\quint{j+1}{q^{-1}}I,&&\rho_j(\pcas)=
\frac{\quint{2j}{q}\quint{2(j+1)}{q^{-1}}}{\quint{2}{q}\quint{2}{q^{-1}}}I,
\end{eqnarray}
but unfortunately he has not been able to come up with a general proof for
either.  Oh well.)

\bigskip

A final comment: if one instead chooses the basis $\{\gen{-},\sqrt{\frac{
\quint{2}{q^{-1}}}{q}}\gen{3},\gen{+}\}$ for the matrix form of the
canonical Killing metric, then
\begin{equation}
\eta_{ab}\propto\left( \begin{array}{ccc}
0&0&\inv{q}\\0&1&0\\q&0&0\end{array}\right) .
\end{equation}
Interestingly, this is the metric for the QG $SO_{q^2}(3)$ \cite{RTF}.
There is some evidence that the classical equivalence between the groups
$SO(3)$ and $SU(2)$ extends to a ``deformed equivalence'' between the QGs
$SO_{q^2}(3)$ and $SU_q(2)$ \cite{SO-SL}, and the result here supports
this.

\section{Conclusions}
\setcounter{equation}{0}

In considering a QHA $\U$, the presence of the universal R-matrix $\R$
implies the existence of the element $u$.  It was then shown in Section
\ref{chap-Killing} that, together with a rep of $\U$, this permits the
introduction of a Killing form with properties (symmetry, invarianve under
the adjoint action) corresponding to those used to define such an object
for the familiar classical case.  If the QHA is also a QLA generated by
$\g$, the Killing metric is just given by plugging two basis elements of
$\g$ into the Killing form.

Under certain assumptions such as invertibility and irreducibility of the
adjoint rep, this Killing metric allows for the definitions of quantities
which are present in the undeformed case, like quadratic casimirs, the
canonical Killing metric, and the index of a rep.  These assumptions are
valid for the classical compact Lie algebras, and it was demonstrated that
they also hold for the particular case of the QLA $\uea{su}{N}$, but that
the classical results are reproduced only for the $(N^2-1)$-dimensional
subspace $\g'$ spanned by the traceless generators.

\bigskip

Therefore, if one wants to construct a physical theory in which the matter
lives in irreps of $\uea{su}{N}$, the states would be labeled by the
quantum numbers corresponding to the values of $\gen{0}$ and $\pcas$ in
that rep, as well as by the weights from the deformed Cartan subalgebra.
For example, in the particular case of $\uea{su}{2}$, the doublet living in
the fundamental irrep (``up'' and ``down'' quarks, maybe?) would consist of
the states (with the ordering $\ket{\gen{0},\pcas,\gen{3}}$)
\begin{equation}
\ket{\pm}:=\ket{-\lambda\quint{\half}{q}\quint{3\over 2}{q^{-1}},
\frac{\quint{3}{q^{-1}}}{\quint{2}{q}\quint{2}{q^{-1}}},\pm\frac{q^{\mp
1}}{q\quint{2}{q^{-1}}}},
\end{equation}
and the triplet living in the adjoint irrep would be
\begin{eqnarray}
\ket{\pm}:=\ket{-\lambda\quint{2}{q^{-1}},\frac{\quint{4}{q^{-1}}}{\quint{2
}{q^{-1}}},\pm q^{\mp 1}},&&\ket{0}:=\ket{-\lambda\quint{2}{q^{-1}},
\frac{\quint{4}{q^{-1}}}{\quint{2}{q^{-1}}},-\lambda}.
\end{eqnarray}

As $q\rightarrow 1$, the latter two quantum numbers approach their
classical values, and the quantum number given by $\gen{0}$ vanishes, as it
must, since for the usual case of {\goth su}$(2)$, $\gen{3}$ and $\pcas$
form a maximal set of commuting operators.  Interestingly, however, there
may exist circumstances where the value of $\gen{0}$ actually removes a
degeneracy; this would occur if, for some value of $q$, there are two
states with the same quantum numbers given by the quadratic casimirs and
weights, but different ones given by $\gen{0}$.  Then even though these two
states would coincide in the $q=1$ case, they are distinct for the deformed
case, and would be distinguishable particles in a theory with a deformed
symmetry group.

\bigskip

In the context of Yang-Mills theories, the formalism here developed will
also come in handy; if $F\equiv\gen{A}F^A$ is the QLA-valued field strength
2-form, and $x=\xi^A\gen{A}\in\g$ (where the coordinates $\xi^A$ are
0-forms), then the fact that $\ad$ is an action implies that
\begin{eqnarray}
\trace{\rho}{\left(x\ad(F\wedge\star F)\right)}&=&\trace{\rho}{
(x_{(1)}\ad F)\wedge\star (x_{(2)}\ad F)}\nonumber\\
&=&\xi^A\epsilon(\gen{A})\trace{\rho}{F\wedge\star F}\nonumber\\
&=&0,
\end{eqnarray}
so if $F\mapsto x\ad F$ describes an infinitesimal gauge transformation, as
it does in the undeformed case, the above trace is gauge-invariant.  Thus,
the quantity $\metric{\rho}{AB}F^A\wedge\star F^B$ may be used in the
action as the kinetic energy of the gauge field.

For the case of $\uea{su}{N}$ (and perhaps others as well), where the
Killing metric is block-diagonal in the primed basis, then by writing
$F=F^0\gen{0}+F^a\pgen{a}$, the action takes the form
\begin{equation}
S=\int\left(\metric{\rho}{ab}F^a\wedge\star F^b+\metric{\rho}{00}F^0\wedge
\star F^0\right)+\ldots,
\end{equation}
which, for the particular example where $N=2$ and $\rho=\adj$, is
\begin{eqnarray}
S&=&\frac{\quint{4}{q^{-1}}}{q^3}\int\left(qF^+\wedge\star F^-+\inv{q}F^-
\wedge\star F^++\frac{q}{\quint{2}{q^{-1}}}F^3\wedge\star F^3\right)
\nonumber\\
&&+\frac{\lambda^2}{q^2}\quint{2}{q^{-1}}^2\quint{3}{q^{-1}}\int F^0\wedge
\star F^0+\ldots .
\end{eqnarray}

Because $\gen{0}$, and thus $\metric{\rho}{00}$, vanishes in the classical
limit, the term involving $F^0$ will not be present in the undeformed case
of {\goth su}$(N)$, but must be included for $q\neq 1$ in order to
guarantee invariance under the adjoint action.  In fact, since the subspace
$\g'$ is the part of $\U$ which survives in the classical limit, the term
quadratic in $F^0$ might be better interpreted as an interaction term in
the lagrangean rather than part of the kinetic energy.  This view is
further supported by the fact that since $\metric{\rho}{00}$ will be small
when $q$ is close to $1$, it may be thought of as a perturbation parameter.

These examples demonstrate that from the point of view of physics, the
existence of a deformed Killing metric (at least for $\uea{su}{N}$)
provides additional structure which may in fact have consequences when
constructing a theory with a quantum symmetry.  In fact, it is possible
that there may be ways to retain some of this extra structure even in the
classical case (\eg if $\lambda F^0$ is kept nonzero in the $\uea{su}{2}$
Yang-Mills example above).  This, however, remains an open question, and
(in the author's opinion) worth further study.

\section*{Acknowledgements}

I would like to thank Chryss Chryssomalakos, Oleg Ogievetsky, Peter Schupp
and Bruno Zumino for helpful comments and suggestions.  I would also like
to mention Olaf Backofen, who brought the connection with \cite{witten} to
my attention.

\newpage

\appendix
\section{Appendix:  The Matrix $D$}
\setcounter{equation}{0}

If $\rho$ is a rep of a QHA $\U$, and $\A$ is the associated QG, the
numerical matrix $D$ is defined to be $\rho(u)$, where $u$ is the element
defined by (\ref{u-def}):
\begin{equation}
\matrix{D}{i}{j}:=\alpha\inprod{u}{\matrix{A}{i}{j}}
\end{equation}
($\alpha$ is just a overall normalization constant).  Several results
follow immediately: first of all, an explicit computation using the
definition of $u$ leads to the result
\begin{equation}
I=\alpha\tr_1(D_1^{-1}\Rhat^{-1})=\alpha^{-1}\tr_2(D_2\Rhat),\label{til-ex}
\end{equation}
where $\tr_J$ is shorthand for the trace over the $J^{th}$ pair of indices,
\eg the $(\matrix{}{i}{j})^{th}$ element of the rightmost expression in the
above equation is $\alpha^{-1}\matrix{D}{m}{n}\matrix{\Rhat}{in}{jm}$.
These relations can be ``inverted'' in the sense of solving them for $D$
and $D^{-1}$ by using the following: given an $N^2\times N^2$ numerical
matrix $M$, one may define the matrix $\tilde{M}=[(M^{t_1})^{-1}]^{t_1}$
($t_J$ denotes transposing with respect to the $J^{th}$ pair of indices)
which satisfies
\begin{equation}
\matrix{M}{im}{n\ell}\matrix{\tilde{M}}{nk}{jm}=\matrix{M}{mi}{\ell n}
\matrix{\tilde{M}}{kn}{mj}=\delta^i_j\delta^k_{\ell}.\label{tilde}
\end{equation}
Hence, by using (\ref{til-ex}), one can find expressions for $D$ and its
inverse:
\begin{eqnarray}
D=\alpha\tr_2(P\tilde{R}),&&D^{-1}=\alpha^{-1}\tr_2(P\widetilde{(R^{-1})}).
\label{D-def}
\end{eqnarray}

Recall that the element $c:=uS(u)$ is central in $\U$; therefore, in an
irrep, $\rho(c)$ must be proportional to the unit matrix.  Hence, define
the constant $\beta$ by means of the identity
\begin{equation}
\inprod{c}{\matrix{A}{i}{j}}=(\alpha\beta)^{-1}\delta^i_j .
\end{equation}
Using the explicit expressions for $c$ and $u$, one obtains
\begin{equation}
I=\beta^{-1}\tr_1(D_1^{-1}\Rhat)=\beta\tr_2(D_2\Rhat^{-1}),
\end{equation}
or, by ``inverting'',
\begin{eqnarray}
D=\beta^{-1}\tr_1(P\widetilde{(R^{-1})}),&&D^{-1}=\beta\tr_1(P\tilde{R}).
\end{eqnarray}

\bigskip

The dual version (in the QG $\A$, that is) of (\ref{uxu}) is
\begin{equation}
S^2(A)=DAD^{-1}.
\end{equation}
This identity, the definition of the $D$-matrix, and (\ref{RAA}) give
\begin{equation}
(D^{-1})^tA^tD^tS(A)^t=S(A)^t(D^{-1})^tA^tD^t =\IA,\label{DADA}
\end{equation}
and (\ref{RAA}) and (\ref{DADA}) together then imply that
\begin{equation}
\tilde{R}=D_1 ^{-1}R^{-1}D_1 =D_2R^{-1}D^{-1}_2 .
\end{equation}
Then from either this relation or the fact that $(S^2\otimes S^2)(\R)=\R$,
it follows that
\begin{equation}
D_1D_2R=RD_1D_2.
\end{equation}

The properties of $D$ just described imply the following: if $M$ is an
$N\times N$ matrix, then
\begin{eqnarray}
\matrix{\tr_1(D_1^{-1}R^{-1}M_1 R)}{i}{j}&=&\matrix{\tr_1(D_1^{-1}R_{21}M_1
R_{21}^{-1})}{i}{j}\nonumber\\
&=&\tr(D^{-1}M)\delta^i_j.
\end{eqnarray}
Also, if the elements of $M$ commute with the elements of $A$,
\begin{equation}
\tr(D^{-1}S(A)MA)=\tr(D^{-1}M).
\end{equation}
In particular, if $M$ is a matrix on which $\A$ right coacts via $\DA
(\matrix{M}{i}{j})=\matrix{M}{k}{\ell}\otimes S(\matrix{A}{i}{k})
\matrix{A}{\ell}{j}$, then (\ref{DADA}) implies
\begin{equation}
\DA(\tr(D^{-1}M))=\tr(D^{-1}M)\otimes\IA.\label{inv-tr}
\end{equation}
For this reason, $\tr(D^{-1}M)$ is called the {\em invariant trace of}
$M$.

\bigskip

For QLAs, this matrix makes its appearance in the following way: the matrix
$\tilde{\real}$, defined by means of (\ref{tilde}) using the numerical
R-matrix $\bigR{}{}$ of the QLA, is given by
\begin{equation}
S(\O{C}{A})\ad\gen{D}=\matrix{\tilde{\real}}{AB}{CD}\gen{B}.\label{Rtilde}
\end{equation}
Notice that (\ref{O-hopf}) implies that
\begin{eqnarray}
S^2(\gen{A})&=&S^2(\O{A}{B})\gen{C}S(\O{B}{C})\nonumber\\
&=&S(\O{A}{B})\ad\gen{B}\nonumber\\
&=&\bigD{B}{A}\gen{B},
\end{eqnarray}
where the numerical matrix $\bigD{}{}$ is defined through (\ref{Rtilde}) as
\begin{equation}
\bigD{A}{B}=\matrix{\tilde{\real}}{CA}{BC},\label{bigD-def}
\end{equation}
which is precisely the same as (\ref{D-def}) in QLA language.  Similarly,
$\bigD{}{}$ satisfies
\begin{eqnarray}
S^2(\O{A}{B})&=&\bigD{C}{A}\O{C}{D}(\field{D}^{-1})^B{}_D,\nonumber\\
\bigD{}{1}\bigD{}{2}\bigR{}{}&=&\bigR{}{}\bigD{}{1}\bigD{}{2},\nonumber\\
\matrix{\tilde{\real}}{AB}{CD}&=&\matrix{(\field{D}^{-1}_1\bigR{-1}{}
\field{D}_2)}{AB}{DC}.
\end{eqnarray}

\end{document}